\documentclass{siamart1116}
\usepackage{amsfonts}
\usepackage{amsmath,amssymb}
\usepackage[all]{xy}
\usepackage{tikz}
\usetikzlibrary{shadings}
\usepackage{enumitem}
\usetikzlibrary{trees}
\usetikzlibrary{arrows,shapes,snakes,automata,backgrounds,petri}
\usepackage{accents}
\usepackage{bm}
\usepackage{nicefrac}

\newcommand{\doublehat}[1]{%
\bar{#1}}

\newcommand{\hatt}[1]{%
\hat{#1}}

\mathchardef\mhyphen="2D 

\newtheorem{thm}{Theorem}
\newtheorem{thm1}{Theorem}[section]
\newtheorem{lemma1}[thm1]{Lemma}

\newtheorem{def1}[thm1]{Definition}
\newtheorem{example}[thm1]{Example}
\newtheorem{notation}[thm1]{Notation}

\theorembodyfont{\normalfont} 
\newtheorem{remark}[thm1]{Remark}
\newtheorem{subproof}[thm1]{Subproof}

\usepackage{pgfplots,comment}

\definecolor{light-gray}{gray}{0.9}
\definecolor{mygray}{gray}{0.8}

\title{Building oscillatory chemical reaction networks by adding reversible reactions}
\author{Murad Banaji\footnotemark[1]}
\begin{document}
\maketitle

\renewcommand{\thefootnote}{\fnsymbol{footnote}}

\footnotetext[1]{Middlesex University, London, Department of Design Engineering and Mathematics. {\tt m.banaji@mdx.ac.uk}.}
\renewcommand{\thefootnote}{\arabic{footnote}}

\begin{abstract}
We show that if a chemical reaction network (CRN) admits nondegenerate (resp., linearly stable) oscillation, and we add new reversible reactions involving new species to this CRN, then the new CRN so created also admits nondegenerate (resp., linearly stable) oscillation provided certain mild and easily checked conditions are met. This claim that the larger CRN ``inherits'' oscillation from the smaller one, provided it is built from the smaller CRN in an appropriate way, follows an analogous result involving multistationarity. It also adds to a number of prior results on the inheritance of oscillation; these collectively often allow us to determine the capacity of a given network for oscillation based on an analysis of its subnetworks. 
\end{abstract}

\section{Introduction and statement of the main result}

A mathematically interesting and practically important question is when we can infer some dynamical behaviour in a network model based on knowledge that this behaviour occurs in some model of a subnetwork. Results in this area focussed on chemical reaction networks (CRNs), for example in \cite{joshishiu,feliuwiufInterface2013,Joshi.2013aa,JoshiShiu2016,banajipanteaMPNE,banajiCRNosci}, have illustrated that this is a rather subtle question. Some intuitively plausible claims turn out to be hard to prove, or to be false. This paper contributes to this literature. The dynamical behaviour of interest here is oscillation, and the main result to be proved (Theorem~\ref{mainthm} below) was conjectured to hold in the concluding sections of \cite{banajiCRNosci}; however the proof turned out to be somewhat harder than expected. 

\textcolor{black}{The work here can be motivated by the following natural questions about sequestration or inhibition. To make the arguments concrete, consider any standard biochemical model exhibiting stable oscillations, such as the classical two-pool calcium model of Goldbeter {\em et al} \cite{GoldbeterTwoPool} or models of glycolytic oscillations \cite{Westermark}. 
Now suppose we introduce a reversible ``sequestration'' process where one of the existing chemical species can bind reversibly to some species, either already in the model or new, to create a new inactive complex. Equivalently (from a mathematical point of view), we introduce an allosteric inhibitor which binds reversibly to some existing species converting it to an inactive form. Could such a process, with certainty, destroy stable oscillation? The answer is no -- maintaining the kinetics of the original reactions, the model will still stably oscillate if we choose mass action kinetics and appropriate rate constants for the inhibition/sequestration process. This claim is intuitively plausible. The idea is that provided (i) the binding and unbinding rates are fast so the new reaction is trying hard to equilibriate, and (ii) the ratio of rates is chosen so the bound, inactive, form has low equilibrium concentration, then the dynamics of the enlarged model projected onto the original species space is ``close'' to that of the original model, and stable oscillation should survive by perturbation arguments.}

\textcolor{black}{Theorem~\ref{mainthm} of this paper is much more general than this motivating discussion suggests, but is guided by the desire to find under what circumstances similar arguments work. The Theorem states that we can build a new oscillatory netork by adding into an existing oscillatory network new reversible reactions, but with a caveat:} some new species must figure nondegenerately in the new reactions. This condition is made precise later, but can easily be illustrated in the special case of a single added reaction, when it becomes: ``{\em there must be a net change in at least one new species in the added reaction}''. Example~\ref{exbasic} provides a simple illustration of the result in this special case. Meanings of the terms, and assumptions about reaction kinetics, will follow later. 

\begin{example}
\label{exbasic}
Suppose that we have a CRN $\mathcal{R}$ on chemical species $X_1, \ldots, X_n$ admitting linearly stable oscillation. We now build a larger CRN $\mathcal{R}'$ by adding to $\mathcal{R}$ a reaction $\mathcal{R}_0$ involving some new species. Then, for example:
\begin{itemize}
\item[(i)] If $\mathcal{R}_0$ is $X_1 + X_{n+1} \rightleftharpoons 2X_{n+1}$ then $\mathcal{R}'$ also admits linearly stable oscillation: there is a net change in the new species $X_{n+1}$ in the added reaction.
\item[(ii)] If $\mathcal{R}_0$ is $X_1+ X_{n+1} \rightleftharpoons X_{n+1} + X_{n+2}$ then $\mathcal{R}'$ admits linearly stable oscillation: there is a net change in the new species $X_{n+2}$ in the added reaction.
\item[(iii)] If $\mathcal{R}_0$ is $X_1 + X_{n+1} \rightleftharpoons X_{n+1}$ then we cannot conclude from Theorem~\ref{mainthm} that $\mathcal{R}'$ admits oscillation: there is a new species involved, but the added reaction does not cause any net change in this new species. 
\end{itemize}
\end{example}

It is little surprise that perturbation theory (both regular and singular) forms the backbone of the proof of Theorem~\ref{mainthm}. The challenge which takes up the majority of our effort here is to recast the basic problem in a form amenable to geometric singular perturbation theory approaches. \textcolor{black}{This requires effort since, firstly, adding reactions into a system alters the dynamics of both the new and the original model species and, secondly, various mathematical necessary conditions for the application of geometric singular perturbation theory are easily violated.}

\textcolor{black}{Returning to the broader context of the work,} this paper contributes to the literature on oscillation in CRNs. This literature has a considerable history and includes theoretical, numerical, and algorithmic work, focussed on both ruling out oscillation, and finding oscillation or bifurcations which give rise to oscillation. It would be hard to compile a complete list of papers about, or with important implications for, oscillations in CRNs. Instead, the following is a small sample, illustrating both numerical and applied work, and some key strands of classical and modern theory: \cite{hornjackson,field1974,feinberg,dicera,goldbeter1995,WolfOsci,gatermann,Qiao.2007aa,minchevaroussel,banajidynsys,domijan,angelileenheersontag,donnellbanaji,abphopf,Woller_2014,errami2015}. Inheritance results of the kind here provide an important theoretical tool for guaranteeing the occurrence of oscillation in CRNs without resorting to numerical investigation. 

Presenting a detailed introduction to the mathematical theory of CRNs can often take up a considerable chunk of a paper, and a minimal approach is adopted here: we intersperse key definitions into the text without much discussion. The reader is referred to some of the papers referenced above and to \cite{banajipantea} for a more thorough background using notation and conventions close to those adopted here. 

We now turn to statement of the main result. Consider a chemical reaction network $\mathcal{R}$ involving $n$ species $X_1, \ldots, X_n$ collectively termed $X$. Let species $X_i$ have concentration $x_i$ ($i=1, \ldots, n$). We are interested in positive concentrations, namely $x:=(x_1, \ldots, x_n)^{\mathrm{t}} \in \mathbb{R}^n_{\gg 0} := \{x\in \mathbb{R}^n \colon x_i > 0\,\,\mbox{for all}\,\, i\}$. Suppose that we have $r_0$ chemical reactions involving $X$, and that the $i$th reaction has reaction vector $\Gamma_i$ and reaction rate $v_i(x)$. We assume only that $v_i\colon \mathbb{R}^n_{\gg 0} \to \mathbb{R}$ is $C^2$; \textcolor{black}{other than this condition, the kinetics of the reactions is arbitrary}. $v(x) := (v_1(x), \ldots, v_{r_0}(x))^{\mathrm{t}}$ is termed the {\em reaction rate vector} for the system, and $\Gamma := [\Gamma_1|\cdots|\Gamma_{r_0}]$ is termed the {\em stoichiometric matrix} of the system. Then the following system of ODEs on $\mathbb{R}^n_{\gg 0}$ describes the evolution of the concentration vector $x$. 
\begin{equation}
\label{eqbasic0}
\dot x = \Gamma v(x)\,,
\end{equation}
Note that the RHS of (\ref{eqbasic0}) belongs to $\mathrm{im}\,\Gamma$, a linear subspace of $\mathbb{R}^n$ termed the {\em stoichiometric subspace} of the system, and consequently cosets of $\mathrm{im}\,\Gamma$ are invariant under the local flow defined by (\ref{eqbasic0}) on $\mathbb{R}^n_{\gg 0}$. The intersection of these cosets of $\mathrm{im}\,\Gamma$ with $\mathbb{R}^n_{\gg 0}$ are termed the {\em positive stoichiometry classes} of the system. 

Now let $m\geq 1$ and $k\geq 0$ be integers and suppose that we add to the system $m$ new reversible reactions involving $m+k$ new species $Y_1, \ldots, Y_{m+k}$, collectively termed $Y$. The new CRN obtained from $\mathcal{R}$ by adding in the new reversible reactions will be termed $\mathcal{R}'$. In order to state a nondegeneracy condition on the added reactions we need to describe these reactions, and for this we introduce some notation.

Given a list of species, say $X_1, \ldots, X_n$, and a nonnegative integer vector of the same length, say $c=(c_1, \ldots, c_n)^{\mathrm{t}}$, we write $c \cdot X$ for the formal sum (i.e., {\em complex} in CRN terminology) $c_1X_1 + c_2 X_2 + \cdots + c_nX_n$. We simply write ``$0$'' for the zero complex \textcolor{black}{$0X_1 + \cdots + 0X_n$}. Using this notation, let the new reactions be:
\begin{equation}
\label{addedreacs}
a_i\cdot X +b_i\cdot Y \rightleftharpoons a_i'\cdot X + b_i'\cdot Y,\quad (i = 1, \ldots, m)\,.
\end{equation}
Here $a_i, a_i', b_i$ and $b_i'$ are nonnegative integer vectors of length $n$, $n$, $m+k$ and $m+k$ respectively. Define $a = (a_1| a_2|\cdots|a_m) \in \mathbb{R}^{n \times m}$, with $a' \in \mathbb{R}^{n \times m}$, $b \in \mathbb{R}^{(m+k) \times m}$ and $b' \in \mathbb{R}^{(m+k) \times m}$ defined similarly. Define $\alpha = a'-a \in \mathbb{R}^{n \times m}$ and $\beta = b'-b \in \mathbb{R}^{(m+k) \times m}$. $\alpha$ records the net stoichiometric changes of the old species $X$ in the added reactions. $\beta$ records the net stoichiometric changes of the new species $Y$ in the added reactions and occurs in a nondegeneracy condition in Theorem~\ref{mainthm} below. Let $y_i$ denote the concentration of $Y_i$ ($i=1, \ldots, m+k$), and define $y:=(y_1, \ldots, y_{m+k})^{\mathrm{t}}$. If the new reactions have reaction rate vector $q\colon \mathbb{R}^n_{\gg 0} \times \mathbb{R}^{m+k}_{\gg 0} \to \mathbb{R}^m$, then $\mathcal{R}'$ evolves according to:

\begin{equation}
\label{perteq}
\left(\begin{array}{c}\dot x\\\dot y\end{array}\right) = \left(\begin{array}{cc}\Gamma&\alpha\\0&\beta\end{array}\right)\left(\begin{array}{c}v(x)\\q(x,y)\end{array}\right)\,
\end{equation}
on $\mathbb{R}^n_{\gg 0} \times \mathbb{R}^{m+k}_{\gg 0}$. We are now ready to state our main result, although some of the definitions to make it precise will follow.

\begin{thm}
\label{mainthm}
Suppose the CRN $\mathcal{R}$ with evolution described by (\ref{eqbasic0}) has a nondegenerate (resp., linearly stable) positive periodic orbit. Let $\mathcal{R}'$ be the CRN with evolution described by (\ref{perteq}), obtained by adding in the reactions of (\ref{addedreacs}) to $\mathcal{R}$. Suppose (i) $\beta$ has rank equal to $m$, its number of columns, and (ii) the added reactions are given mass action kinetics. Then rate constants can be chosen for the added reactions in such a way that $\mathcal{R}'$ has a nondegenerate (resp., linearly stable) positive periodic orbit. 
\end{thm}

By a {\em positive} periodic orbit, we mean a periodic orbit that lies in the (strictly) positive orthant. By a {\em nondegenerate} periodic orbit, we mean one that is hyperbolic relative to its stoichiometry class. {\em Linearly stable} is also taken to mean linearly stable relative to its stoichiometry class. Precise statement of these latter conditions is deferred to Section~\ref{sectech}.

Theorem~\ref{mainthm} is exactly analogous (including the condition that $\beta$ has rank $m$) to Theorem~5 in \cite{banajipanteaMPNE}, replacing ``multiple positive nondegenerate (resp., linearly stable) equilibria'' with ``a nondegenerate (resp., linearly stable) positive periodic orbit''. The proof draws heavily both on techniques developed in \cite{banajipanteaMPNE}, and on singular perturbation theory approaches which formed the basis for the proof of Theorem~4 in \cite{banajiCRNosci}. It is hoped that the proof of Theorem~\ref{mainthm} provides a template for, \textcolor{black}{or at least insight towards,} the proof of further inheritance results on CRNs.

\section{An example}
\label{secexample}

The proof of Theorem~\ref{mainthm} is constructive: it not only tells us about inheritance of oscillation, but also gives information about parameter regions at which oscillation occurs. Before the proof, we present an example illustrating the result, including how to choose parameter values at which we can observe inherited oscillation. \textcolor{black}{We will re-examine in detail this example in Section~\ref{seceluc} after presenting the proof of Theorem~\ref{mainthm} and use it to help elucidate various ideas and quantities in the proof. Here we use it only to outline the conclusions.}

\begin{figure}[h]
\begin{minipage}{0.48\textwidth}
\begin{center}
\includegraphics[scale=.3]{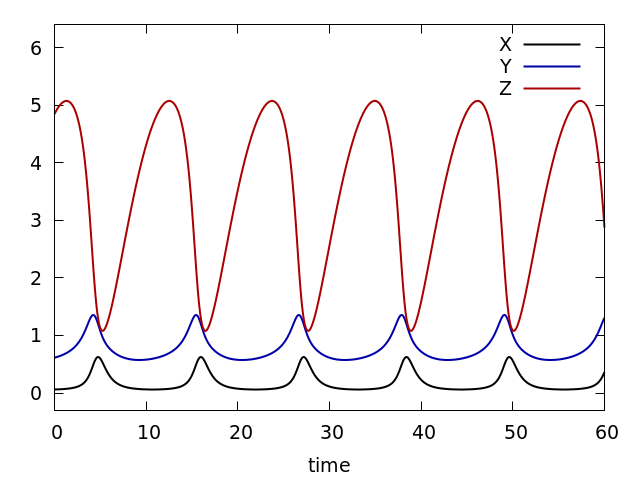}
\end{center}
\end{minipage}
\hfill
\begin{minipage}{0.48\textwidth}
\begin{center}
\includegraphics[scale=.3]{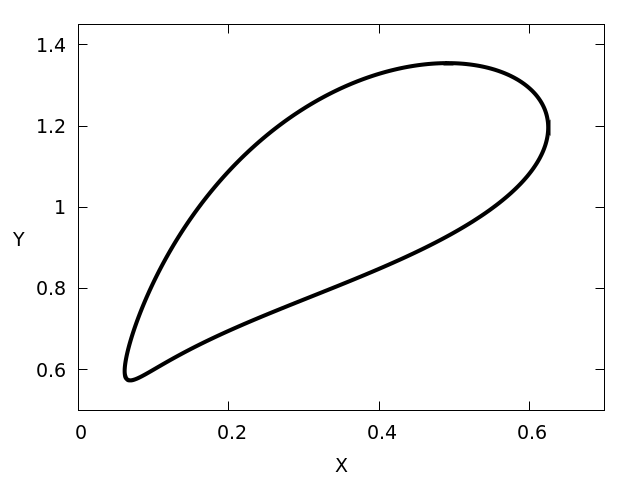}
\end{center}
\end{minipage}
\caption{\label{fig1} Simulation of the CRN $\mathcal{R}_1$ with mass action kinetics and rate constants as given in the text. {\em Left.} Evolution of the concentrations of $X$, $Y$ and $Z$. {\em Right.} The projection of the periodic orbit onto $X\mhyphen Y$ coordinates.}
\end{figure}

In \cite{banajiCRNosci}, the following was presented as an example of a CRN which admits stable oscillation with mass action kinetics:
\begin{equation}
X+Z \overset{\scriptstyle{k_1}}\longrightarrow 2Y \overset{\scriptstyle{k_2}}\longrightarrow X+Y, \quad 0 \overset{\scriptstyle{k_3}}{\underset{\scriptstyle{k_4}}\rightleftharpoons} X, \quad 0 \overset{\scriptstyle{k_5}}{\underset{\scriptstyle{k_6}}\rightleftharpoons} Y, \quad 0 \overset{\scriptstyle{k_7}}{\underset{\scriptstyle{k_8}}\rightleftharpoons} Z\,. \tag{\mbox{$\mathcal{R}_1$}}
\end{equation}
This is an example of a so-called {\em fully open} CRN on three species $X, Y$ and $Z$, as it includes the inflow-outflow reactions $0 \rightleftharpoons X$, $0 \rightleftharpoons Y$ and $0 \rightleftharpoons Z$. The reactions are labelled with their rate constants. In numerical simulations we easily find parameter regions where the CRN admits a stable periodic orbit. For example, setting $k_1=4,\, k_2=3,\, k_3=0.2,\, k_4=2,\, k_5=0.3,\, k_6=2.5,\, k_7=2.5$ and $k_8=0.2$, and choosing initial conditions $X_0=Y_0=Z_0=1$ we find the system settles, after initial transients, onto the periodic orbit shown in Figure~\ref{fig1}. We assume that the system does indeed have a positive, linearly stable periodic orbit at these parameter values.

Now suppose that we add in two further reversible reactions $Y \rightleftharpoons U+V$ and $U+X \rightleftharpoons 2V+W$ involving three new species $U, V$ and $W$ to obtain the system
\begin{equation}
X+Z \overset{\scriptstyle{k_1}}\longrightarrow 2Y \overset{\scriptstyle{k_2}}\longrightarrow X+Y, \,\,\, 0 \overset{\scriptstyle{k_3}}{\underset{\scriptstyle{k_4}}\rightleftharpoons} X, \,\,\, 0 \overset{\scriptstyle{k_5}}{\underset{\scriptstyle{k_6}}\rightleftharpoons} Y, \,\,\, 0 \overset{\scriptstyle{k_7}}{\underset{\scriptstyle{k_8}}\rightleftharpoons} Z, \,\,\, Y \overset{\scriptstyle{k_9}}{\underset{\scriptstyle{k_{10}}}\rightleftharpoons} U+V, \,\,\, U+X \overset{\scriptstyle{k_{11}}}{\underset{\scriptstyle{k_{12}}}\rightleftharpoons} 2V+W\,. \tag{\mbox{$\mathcal{R}_2$}}
\end{equation}
The matrix $\beta$ representing the net stoichiometric changes of the new species in the added reactions is, in this case 
\[
\left(\begin{array}{rr}1&-1\\1&2\\0&1\end{array}\right)
\]
which clearly has rank $2$. Consequently, the nondegeneracy condition in Theorem~\ref{mainthm} is satisfied, and the theorem tells us that $\mathcal{R}_2$ admits stable oscillation with mass action kinetics.

\begin{figure}[h]
\begin{minipage}{0.48\textwidth}
\begin{center}
\includegraphics[scale=.3]{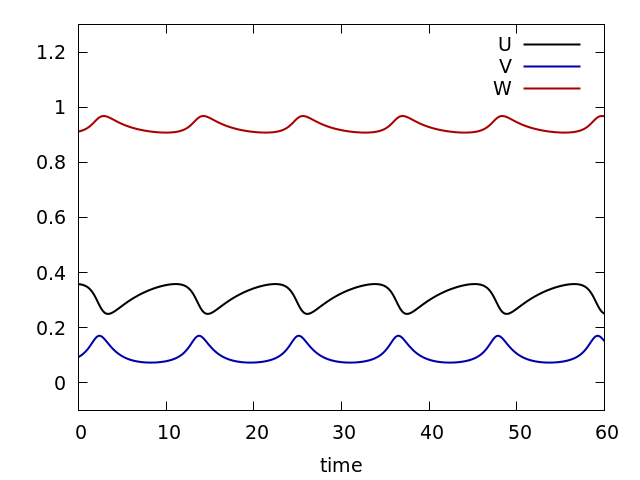}
\end{center}
\end{minipage}
\hfill
\begin{minipage}{0.48\textwidth}
\begin{center}
\includegraphics[scale=.3]{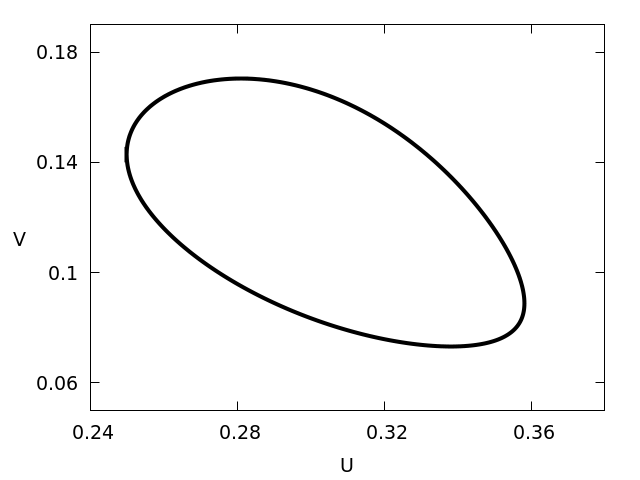}
\end{center}
\end{minipage}
\caption{\label{fig2} Simulation of the CRN $\mathcal{R}_2$ with mass action kinetics and rate constants as given in the text. {\em Left.} Evolution of the concentrations of $U$, $V$ and $W$. {\em Right.} The projection of the periodic orbit onto $U\mhyphen V$ coordinates.}
\end{figure}

The proof of Theorem~\ref{mainthm} also tells us \textcolor{black}{how to set rate constants for the added reactions to obtain} oscillation in $\mathcal{R}_2$. We define two parameters $\epsilon$ and $\eta$, set $k_9=\epsilon^{-1}$, $k_{10}=\epsilon^{-1}\eta^{-2}$, $k_{11}=\epsilon^{-1}\eta^{-1}$ and $k_{12}=\epsilon^{-1}\eta^{-2}$, and leave all other rate constants as before; then the proof of Theorem~\ref{mainthm} tells us that by choosing and fixing $\eta>0$ sufficiently small, and subsequently choosing and fixing $\epsilon>0$ sufficiently small, we can ensure that $\mathcal{R}_2$ has a positive periodic orbit which is linearly stable relative to its stoichiometry class. Moreover, with these choices, variation in the values of $U$, $V$ and $W$ on the periodic orbit will be small; the values of $U$ and $V$ on the periodic orbit will be small; the values of $W$ on the periodic orbit can be controlled by the choice of initial data; and the values of $X$, $Y$ and $Z$ on the periodic orbit will be close to their original values in the absence of the added reactions. Some plots of the periodic orbit (omitting transient behaviour) in the case $\epsilon=\eta=0.2$, and with initial values of the new variables $U_0=V_0=0$, $W_0=1$ are shown in Figure~\ref{fig2}. Note that $\mathcal{R}_2$ now has a conserved quantity $3W+U-V$. 

\section{Technical preliminaries}
\label{sectech}
Before presenting the proof of Theorem~\ref{mainthm} we need some notational and mathematical preliminaries from analysis and the theory of differential equations.

\subsection{Basic notation, conventions, and definitions} 
We draw heavily on \cite{banajipanteaMPNE} and \cite{banajiCRNosci} here. 

\begin{def1}[Positive subsets of $\mathbb{R}^n$]
We refer to a subset of $\mathbb{R}^n$ as positive if it is a subset of 
\[
\mathbb{R}^n_{\gg 0}:=\{x \in \mathbb{R}^n\colon x_i > 0\,\, (i=1, \ldots, n)\}\,.
\]
Given $x \in \mathbb{R}^n$, we write $x \gg 0$ to mean $x \in \mathbb{R}^n_{\gg 0}$. We also define
\[
\mathbb{R}^n_{\geq 0}:=\{x \in \mathbb{R}^n\colon x_i \geq 0\,\, (i=1, \ldots, n)\}\,.
\]
\end{def1}

\begin{notation}[Vector of ones, identity matrix]
$\mathbf{1}$ denotes a vector of ones whose length is inferred from the context. If $\eta$ is any real constant, then $\bm{\eta}$ denotes a vector whose entries are all $\eta$ and whose length is inferred from the context. $I_n$ refers to the $n \times n$ identity matrix.
\end{notation} 

\begin{def1}[Empty vectors and matrices]
In order to simplify some arguments, we formally allow vectors and matrices to be {\em empty}. An empty matrix is one with zero rows, zero columns, or both: when we define a matrix to be $n \times m$, one or both of $n$ or $m$ is allowed to be zero. An empty vector, for example, is a $0 \times 1$ matrix. Empty vectors and matrices obey the following natural rules. (i) If $A$ is an $n \times m$ matrix and $B$ is an $m \times k$ matrix, then $AB$ is defined and is an $n \times k$ matrix, even if some of $n,m$ or $k$ are zero. If $m=0$, but $n$ and $k$ are nonzero, then $AB$ is defined to be the $n \times k$ zero matrix. (ii) Any equation, inequality, or claim involving empty vectors or matrices is vacuously satisfied (provided that it makes sense, dimensionally). (iii) Given an empty vector $y$ and a $k \times 0$ matrix $A$, $y^A$ is defined to be $\mathbf{1}$, a vector of ones of length $k$. 
\end{def1}

\begin{notation}[Sum of a point and a set]
Given a point $x_0\in \mathbb{R}^n$ and a set $A \subseteq \mathbb{R}^n$, $x_0+A$ means, naturally, the following subset of $\mathbb{R}^n$: $\{x_0+y\colon y \in A\}$.
\end{notation}

\begin{notation}[Open ball in $\mathbb{R}^n$]
For any $x\in \mathbb{R}^n$ and $r > 0$, let $B_r(x)$, be the open ball in $\mathbb{R}^n$ of radius $r$ with centre $x$, namely $B_r(x):= \{y \in \mathbb{R}^n \colon |y-x| < r\}$. If the argument $x$ is omitted, it is taken to be zero. The dimension $n$ is to be inferred from the context. 
\end{notation}

\begin{notation}[Hausdorff distance]
Given two nonempty sets $A$ and $B$ in $\mathbb{R}^n$ with the Euclidean metric, $d_{\mathrm{H}}(A, B)$ denotes the Hausdorff distance between $A$ and $B$.
\end{notation}

\begin{notation}[Monomials, vector of monomials]
\label{notmon}
Given $x=(x_1,\ldots, x_n)^{\mathrm{t}}$ and $a = (a_1,\ldots, a_n)$, $x^a$ is an abbreviation for the (generalised) monomial $\prod_ix_i^{a_i}$. If $A$ is an $m \times n$ matrix with rows $A_1, \ldots, A_m$, then $x^A$ means the vector of (generalised) monomials $(x^{A_1}, x^{A_2}, \ldots, x^{A_m})^{\mathrm{t}}$. 
\end{notation} 

\begin{notation}[Entrywise product and entrywise functions]
Given two matrices $A$ and $B$ with the same dimensions, $A \circ B$ will refer to the entrywise (or Hadamard) product of $A$ and $B$, namely $(A\circ B)_{ij} = A_{ij}B_{ij}$. When we apply functions such as $\ln(\cdot)$ and $\exp(\cdot)$ with a vector or matrix argument, we mean entrywise application. Similarly, if $x = (x_1, \ldots, x_n)^{\mathrm{t}}$ and $y = (y_1, \ldots, y_n)^{\mathrm{t}}$, then $x/y$ means $(x_1/y_1, x_2/y_2,\ldots, x_n/y_n)^{\mathrm{t}}$.
\end{notation}

\begin{def1}[Mass action kinetics, rate constants]
A chemical reaction $a\cdot X \rightarrow b\cdot X$ is said to have mass action kinetics if the rate of reaction is $k x^{a^{\mathrm{t}}}$ for some positive constant $k$ termed the {\em rate constant} of the chemical reaction. 
\end{def1}

\begin{notation}[Preimages of sets]
Given a function $f \colon A \to B$, and any $\mathcal{B} \subseteq B$, $f^{-1}(\mathcal{B})$ refers, naturally, to $\{a\in A \colon f(a) \in \mathcal{B}\}$.
\end{notation}

\begin{remark}[Differentiability of functions]
\label{remdifffunc}
When we refer to a function $f$ as being $C^r$ on some set $U \subseteq \mathbb{R}^n$, not necessarily open, we mean that there exists a function $\hatt{f}$ defined and $C^r$ on an open set $V \subseteq \mathbb{R}^n$ containing $U$ and such that $\hatt{f}$ coincides with $f$ on $U$.
\end{remark}

\begin{notation}[Derivatives of functions]
Given a differentiable function $f\colon U\subseteq \mathbb{R}^n \to \mathbb{R}^m$, $Df$ refers both to the derivative of $f$ and also its matrix representation where the bases on $\mathbb{R}^n$ and $\mathbb{R}^m$ are the standard bases or are to be inferred from the context. Given a set of positive integers $n_1, \ldots, n_m, k$, variables $x_i \in \mathbb{R}^{n_i}$ ($i=1, \ldots, m$) and a differentiable function $f \colon \mathbb{R}^{n_1} \times \cdots \times \mathbb{R}^{n_m} \to \mathbb{R}^k$, $D_{x_i}f$ refers to the derivative of $f$ w.r.t. the variable $x_i$ and also its matrix representation. We may also write $D_if$ for the derivative of a function $f$ w.r.t. its $i$th argument, or the matrix representation of this derivative.
\end{notation}

The following three examples, reproduced or adapted from  \cite{banajiCRNosci}, demonstrate how entrywise and monomial notation greatly abbreviate otherwise lengthy calculations.
\begin{example}[Rules of exponentiation]
Let $x,y \in \mathbb{R}^m_{\gg 0}$, $A,B \in \mathbb{R}^{n \times m}$ and $C \in \mathbb{R}^{k \times n}$. Let $O$ refer to the $n \times m$ matrix of zeros. Then (i) $x^O = \mathbf{1}$, (ii) $x^{A+B} = x^A \circ x^B$, (iii) $x^A\circ y^A = (x\circ y)^A$ and (iv) $(x^A)^C = x^{CA}$.
\end{example}

\begin{example}[Logarithm of monomials]
Suppose $x \in \mathbb{R}^m_{\gg 0}$, $y_i \in \mathbb{R}^{n_i}_{\gg 0}$ ($i=1, \ldots, k$) and $A_i \in \mathbb{R}^{m \times n_i}$ ($i=1, \ldots, k$). If $w=x \circ y_1^{A_1} \circ \cdots \circ y_k^{A_k}$, then $\ln w = \ln x + A_1\ln y_1 + \cdots + A_k \ln y_k$. 
\end{example}

\begin{example}[Differentiation of monomials]
\label{exdiffmon}
Suppose $k \in \mathbb{R}^m_{\gg 0}$, $x \in \mathbb{R}^n_{\gg 0}$, $A \in \mathbb{R}^{m \times n}$. Let $w \colon \mathbb{R}^n_{\gg 0} \to \mathbb{R}^m_{\gg 0}$ be defined by $w(x) := k \circ x^A$. Then $Dw = \mathrm{diag}(w)\,A\,\mathrm{diag}(\mathbf{1}/x)$. 
\end{example}

\subsection{Periodic orbits}

We need a number of standard results on periodic orbits and Floquet theory largely as described in Section~2 of \cite{banajiCRNosci}. We summarise these here, but the reader is referred to \cite{banajiCRNosci} and the original sources (\cite{HaleOsci} for example) for more detail. 

Consider some system of ODEs $\dot x = F(x)$ on an open set $U \subseteq \mathbb{R}^n$, satisfying conditions for existence and uniqueness of solutions and hence defining a local flow $\Phi_t\colon U \to U$. For each $x \in U$, $t$ belongs to an open interval including $0$ which in general depends on $x$, and $\Phi_t(x)$ is the point that initial condition $x$ ``reaches'' at time $t$. Given some $T>0$ the orbit of a nontrivial $T$-periodic solution of the ODE system is termed a {\em periodic orbit} of the system. Associated with any such periodic orbit are its {\em Floquet multipliers} (or {\em characteristic multipliers}), namely eigenvalues of $D\Phi_T(x_0)$ where $x_0$ is any point on the periodic orbit, and  $D\Phi_T(x_0)$ refers to the derivative of $\Phi_T$ w.r.t. $x$ evaluated at $x_0$. Here $D\Phi_t(x_0)$ can be regarded as the fundamental matrix solution of the $T$-periodic variational equation $\dot y = DF(\Phi_t(x_0))y$ satisfying $y_0 = I_n$. The choice of $x_0$ does not affect the Floquet multipliers. 

Any periodic orbit always has one Floquet multiplier of $1$ corresponding to the direction tangential to the periodic orbit; the remaining Floquet multipliers are termed the {\em nontrivial} Floquet multipliers of the periodic orbit. If none of the nontrivial Floquet multipliers lie on the unit circle in the complex plane, then the periodic orbit is {\em hyperbolic}, and, in our terminology here, {\em nondegenerate}. If, further, all of the nontrivial Floquet multipliers lie inside the unit circle in the complex plane, then the periodic orbit is {\em linearly stable} and attracts a neighbourhood of itself.

Given any nondegenerate (resp., linearly stable) periodic orbit, by regular perturbation theory arguments involving, for example, the construction of a Poincar\'e map, a nearby nondegenerate (resp., linearly stable) periodic orbit (with nearby period) exists for all ODEs on $U$ $C^1$-close to $\dot x = F(x)$. More precisely, we have the following result, which appears as Lemma~2.1 in \cite{banajiCRNosci}. A proof can be found in Section IV of \cite{Fenichel79}. 
\begin{lemma1}
\label{lemreg}
Let $U \subseteq \mathbb{R}^r$ be open, $\epsilon' > 0$ and $F\colon U \times (-\epsilon', \epsilon')\to \mathbb{R}^r$ be $C^1$. Consider the $\epsilon$-dependent family of ODEs on $U$
\begin{equation}
\label{eqnFloqe}
\dot x = F(x,\epsilon)\,.
\end{equation} 
Suppose that $\dot x = F(x,0)$ has a nontrivial hyperbolic (resp., linearly stable) $T$-periodic orbit $\mathcal{O} \subseteq U$. Then there exists $\epsilon_0 \in (0, \epsilon']$ s.t. for $\epsilon \in (-\epsilon_0, \epsilon_0)$ (\ref{eqnFloqe}) has a hyperbolic (resp., linearly stable) periodic orbit $\mathcal{O}_\epsilon$ satisfying $\lim_{\epsilon \to 0} d_{\mathrm{H}}(\mathcal{O}_\epsilon, \mathcal{O}) = 0$ and with period $T_\epsilon$ satisfying $\lim_{\epsilon \to 0}T_\epsilon = T$. 
\end{lemma1}

An important basic observation that we will frequently need is that the Floquet multipliers of a periodic orbit are invariant under $C^1$-diffeomorphisms, and hence the notions of ``nondegeneracy'' and ``linear stability'' of a periodic orbit are invariant under $C^1$-diffeomorphisms. To be more precise: 
\begin{lemma1}[Invariance of Floquet multipliers]
\label{lemfloqdiff}
Suppose we have a $C^1$ local flow $\Phi_t$ on an open set $U \subseteq \mathbb{R}^n$ with a periodic orbit $\mathcal{O} \subseteq U$. Let $h \colon U \to h(U):=V \subseteq \mathbb{R}^n$ be a $C^1$ diffeomorphism, so that we have a new local flow $\Psi_t = h\circ\Phi_t\circ h^{-1}$ on $V$, and let $\mathcal{O}':=h(\mathcal{O})$ be the corresponding periodic orbit of $\Psi_t$. Then $\mathcal{O}$ and $\mathcal{O}'$ have the same Floquet multipliers. 
\end{lemma1}
\begin{proof}
Given any point $x_0 \in \mathcal{O}$ and the corresponding point $y_0 = h(x_0) \in \mathcal{O}'$, $D\Phi_t(x_0)$ and $D\Psi_t(y_0)$ are linear operators from $T_{x_0}U \cong \mathbb{R}^n$ to $T_{\Phi_t(x_0)}U \cong \mathbb{R}^n$, and $T_{y_0}V \cong \mathbb{R}^n$ to $T_{\Psi_t(y_0)}V \cong \mathbb{R}^n$ respectively. They are defined for all $t \in \mathbb{R}$ (since $\Phi_t(x_0)$ is periodic in $t$), and clearly satisfy $D\Phi_0=\mathrm{id}$ and $D\Psi_0=\mathrm{id}$. Applying the chain rule to $\Psi_t = h\circ\Phi_t\circ h^{-1}$ gives
\[
D\Psi_t(h(x_0)) = Dh(\Phi_t(x_0))\,\,D\Phi_t(x_0)\,\,[Dh^{-1}(h(x_0))]\,.
\]
Setting $t=T$, and observing that $h(x_0)=y$, $\Phi_T(x_0)=x_0$ and $Dh^{-1}(h(x_0)) = [Dh(x_0)]^{-1}$ gives:
\[
D\Psi_T(y_0) = [Dh(x_0)]\,\,D\Phi_T(x_0)\,\,[Dh(x_0)]^{-1}\,.
\]
Now the Floquet multipliers of $\mathcal{O}$ are the eigenvalues of $D\Phi_T(x_0)$ and the Floquet multipliers of $\mathcal{O}'$ are the eigenvalues of $D\Psi_T(y_0)$. Since the final equation shows that  $D\Phi_T(x_0)$ and $D\Psi_T(y_0)$ are similar, the two sets of Floquet multipliers are equal. 
\end{proof}

\begin{remark}[Floquet multipliers relative to a given set]
\label{remfloqrel}
Let $\Phi_t$ be a local flow on an open set $U' \subseteq \mathbb{R}^n$, let $U\subseteq U'$ be locally invariant under $\Phi_t$, and let $\mathcal{O}\subseteq U$ be a periodic orbit of $\Phi_t$. Let $V \subseteq \mathbb{R}^m$ be open and suppose that $h\colon U\to V$ is a $C^1$-diffeomorphism. In the light of Lemma~\ref{lemfloqdiff} it makes sense to refer to the Floquet multipliers of $\mathcal{O}$ relative to $U$. We mean the Floquet multipliers of $h(\mathcal{O})$ for the derived flow $\Psi_t=h\circ\Phi_t \circ h^{-1}$ on $V$, which, by Lemma~\ref{lemfloqdiff}, do not depend on $V$ or $h$.
\end{remark}

{\bf Nondegenerate/linearly stable periodic orbits for a CRN.} Suppose now that we have a chemical reaction network with stoichiometric matrix $\Gamma$ defining a system of ODEs $\dot x = \Gamma v(x)$ as in (\ref{eqbasic0}). Since cosets of $\mathrm{im}\,\Gamma$ are invariant under the local flow defined by such a system, any periodic orbit must belong to one of these sets. If $\Gamma$ has rank $r$, less than its number of rows $n$, then it is easily seen that no periodic orbit can be nondegenerate or linearly stable in the sense described above since any periodic orbit has $n-r$ nontrivial Floquet multipliers with value $1$ corresponding to directions transverse to the coset of $\mathrm{im}\,\Gamma$ on which it lies. In this situation, we follow \cite{banajiCRNosci} and overload the terms {\em nondegenerate} and {\em linearly stable} as follows. We say that a periodic orbit $\mathcal{O}$ is nondegenerate (resp., linearly stable) if it has $r-1$ Floquet multipliers which are disjoint from (resp., inside) the unit circle, or equivalently if it is nondegenerate (resp., linearly stable) {\em relative} to the coset of $\mathrm{im}\,\Gamma$ on which it lies in the sense of Remark~\ref{remfloqrel}. This abuse of terminology should cause no confusion.

\subsection{Some results from analysis}

We need the following form of the implicit function theorem (IFT):
\begin{lemma1}[Implicit Function Theorem]
\label{lemIFT}
Let $W \subseteq \mathbb{R}^n \times \mathbb{R}^m$ be open and $F\colon W \to \mathbb{R}^m$ be $C^r$ ($r \geq 1$). Suppose $F(a,b)=0$ for some $(a,b) \in W$ and the Jacobian matrix $D_2F(a,b)$ (namely with respect to the second variables) is nonsingular. Then there exist $U \subseteq \mathbb{R}^n$, $V \subseteq \mathbb{R}^m$ both open, with $(a,b) \in U\times V \subseteq W$, and a $C^r$ function $\phi \colon U \to V$ satisfying $\phi(a) = b$, and such that
\[
\{(x,y) \in U \times V\,\colon\, F(x,y)=0\} = \{(x,\phi(x))\,\colon\, x \in U\}\,,
\]
namely, the zero set of $F$ in $U \times V$ is precisely the graph of $\phi$.
\end{lemma1}
\begin{proof}
See, for example, Chapter~5 of \cite{pugh2003real}.
\end{proof}

The reader may easily verify that the sets $U$ and $V$ in the statement of the IFT may, without loss of generality, be chosen to be open balls with centres $a$ and $b$ in $\mathbb{R}^n$ and $\mathbb{R}^m$ respectively. 

We need the following consequence of the IFT. 

\begin{lemma1}[IFT extended to a compact set]
\label{lemIFT1}
Let $W \subseteq \mathbb{R}^n \times \mathbb{R}^m$ be open and $F\colon W \to \mathbb{R}^m$ be $C^r$ ($r \geq 1$). Let $X$ be a compact set in $\mathbb{R}^n$ such that $X \times \{0\} \subseteq W$. Suppose that, for all $x \in X$, $F(x,0)=0$ and the Jacobian matrix $D_2F(x,0)$ is nonsingular for each $x \in X$. Then there exist an open set $U \subseteq \mathbb{R}^n$ containing $X$, $t > 0$ such that $U \times B_t \subseteq W$, and a $C^r$ function $\phi\colon U \to B_t$ whose graph is precisely the zero-set of $F$ in $U \times B_t$, namely,
\[
\{(x,y) \in U \times B_t \,\colon\, F(x,y)=0\} = \{(x,\phi(x))\,\colon\, x \in U\}\,.
\]
\end{lemma1}
\begin{proof}
We apply the IFT at $(x,0)$ for each $x \in X$. For each $x \in X$, there exist $s_x>0$ and $t_x>0$ such that $B_{s_x}(x) \times B_{t_x} \subseteq W$, and a $C^r$ function $\phi_x\colon B_{s_x}(x) \to B_{t_x}$ such that the zero set of $F$ in $B_{s_x}(x) \times B_{t_x}$ is precisely the graph of $\phi_x$, namely,
\[
\{(x,y) \in B_{s_x}(x) \times B_{t_x}\,\colon\, F(x,y)=0\} = \{(x,\phi_x(x))\,\colon\, x \in B_{s_x}(x)\}\,.
\]
We now choose a finite set $\{x_i\} \subseteq X$ such that $U':=\cup B_{s_{x_i}}(x_i)$ forms an open cover of $X$. We define the function $\hatt{\phi}\colon U' \to \mathbb{R}^m$ via $\hatt{\phi}(x) = \phi_{x_i}(x)$ where $x_i$ is chosen as {\em any} element such that $x \in B_{s_{x_i}}(x_i)$. $\hatt{\phi}$ is a well defined function since if $x \in B_{s_{x_i}}(x_i) \cap B_{s_{x_j}}(x_j)$, then $\phi_{x_i}(x) = \phi_{x_j}(x)$ (as $\phi_{x_i}(x)$ and $\phi_{x_j}(x)$ must certainly both lie in one of $B_{t_{x_i}}$ or $B_{t_{x_j}}$). It is also clear that $\hatt{\phi}$ is $C^r$ since it coincides with the $C^r$ functions $\{\phi_{x_i}\}$. Let $t := \min\{t_{x_i}\}$. Since $\hatt{\phi}$ is continuous, $X$ is compact, and $\hatt{\phi}(x)=0$ for $x \in X$, there exists an open neighbourhood $U\subseteq U'$ of $X$ such that $x \in U$ implies $|\hatt{\phi}(x)| < t$. Define $\phi:= \left.\hatt{\phi}\right|_{U}$. Clearly $\phi$ satisfies the claims of the lemma, and in particular $\{(x,y) \in U \times B_t\,\colon\, F(x,y)=0\} = \{(x,\phi(x))\,|\, x \in U\}\,.$
\end{proof}

We will need the following technical lemma in order to make uniform estimates on compact sets. Notation is fixed to be consistent with those proofs where Lemma~\ref{lemuniform} is used.
\begin{lemma1}
\label{lemuniform}
Let $\mathcal{Z} \subseteq \mathbb{R}^r$ be compact, $m$ a positive integer, and $\eta'>0$ a positive constant. Let $\theta \colon \mathcal{Z} \times [-\eta',\eta'] \to \mathbb{R}^m$ satisfy
\begin{itemize}
\item $\theta$ is $C^1$ with Lipschitz continuous derivative on its domain of definition (see Remark~\ref{remdifffunc}). For example, if $\theta$ is $C^2$, then this condition certainly holds.
\item $\theta(z, 0)=0$ for all $z \in \mathcal{Z}$. 
\end{itemize}
Then $\theta/\eta\colon \mathcal{Z} \times [-\eta',\eta']\backslash\{0\} \to \mathbb{R}^m$ has a continuous extension to $\mathcal{Z} \times [-\eta', \eta']$. \textcolor{black}{Explicitly}, $\hatt{\theta} \colon \mathcal{Z} \times [-\eta',\eta'] \to \mathbb{R}^m$ defined by
\[
\hatt{\theta}(z,\eta):=\left\{\begin{array}{rcl}\theta(z,\eta)/\eta & (\eta\neq 0)\\D_\eta \theta(z, 0) & (\eta=0)\end{array}\right.
\]
is continuous on $\mathcal{Z} \times [-\eta',\eta']$. Consequently $R(z, \eta):=\hatt{\theta}(z,\eta) - D_\eta \theta(z, 0)$ is continuous on $\mathcal{Z} \times [-\eta',\eta']$. 
\end{lemma1}
\begin{proof}
It is trivial that $\hatt{\theta}$ is continuous (in fact, $C^1$) at points in its domain where $\eta \neq 0$. So we need to show that it is continuous at an arbitrary point of $\mathcal{Z} \times \{0\}$. Define $R(z, \eta):=\hatt{\theta}(z,\eta) - D_\eta \theta(z, 0)$ on $\mathcal{Z} \times [-\eta',\eta']$. Then $R(z,0)=0$ and $\theta(z, \eta) = \eta[D_\eta \theta(z, 0) + R(z, \eta)]$ is an identity; on the other hand, Taylor's theorem tells us that, \textcolor{black}{for any fixed $z$}, $\lim_{\eta \to 0} R(z, \eta)=0$. We would like to show that given any $z_0 \in \mathcal{Z}$, $|\hatt{\theta}(z,\eta) - \hatt{\theta}(z_0, 0)| \to 0$ as $(z,\eta) \to (z_0,0)$.  If $\eta=0$, then
\[
|\hatt{\theta}(z,\eta) - \hatt{\theta}(z_0, 0)| = |\hatt{\theta}(z,0) - \hatt{\theta}(z_0, 0)| = |D_\eta \theta(z, 0) - D_\eta \theta(z_0, 0)|
\]
which, by continuity of $D_\eta\theta$ at $(z_0,0)$ can be made arbitrarily small by choosing $|z-z_0|$ sufficiently small. So now consider the case $\eta \neq 0$. We assume that $\eta>0$; the case $\eta < 0$ requires minor modifications below. By the triangle inequality:
\begin{eqnarray}
\nonumber
|\hatt{\theta}(z,\eta) - \hatt{\theta}(z_0, 0)| & \leq &|\hatt{\theta}(z,\eta) - \hatt{\theta}(z_0, \eta)| + |\hatt{\theta}(z_0, \eta)-\hatt{\theta}(z_0,0)|\\
\label{eqprf0}
&=&\frac{1}{\eta}|\theta(z,\eta) - \theta(z_0, \eta)| + |\hatt{\theta}(z_0, \eta)-\hatt{\theta}(z_0,0)|\,.
\end{eqnarray}
The final term is simply the magnitude of the remainder $R(z_0, \eta)$ in the Taylor expansion:
\begin{equation}
\label{eqprf1}
|\hatt{\theta}(z_0, \eta)-\hatt{\theta}(z_0,0)| = |\theta(z_0,\eta)/\eta-D_\eta \theta(z_0, 0)| = |R(z_0, \eta)|\,.
\end{equation}
On the other hand, using the fundamental theorem of calculus, we get
\[
|\theta(z,\eta) - \theta(z_0, \eta)| = \left|\int_0^\eta (D_\eta \theta(z, \sigma) - D_\eta \theta(z_0, \sigma))\,\mathrm{d}\sigma\right| \leq \int_0^{\eta} \left|D_\eta \theta(z, \sigma) - D_\eta \theta(z_0, \sigma)\right|\,\mathrm{d}\sigma\,.
\]
Lipschitz continuity of $D_\eta \theta$ means that $\left|(D_\eta \theta(z, \sigma) - D_\eta \theta(z_0, \sigma))\right| \leq K|z-z_0|$, where $K$ is the Lipschitz constant of $D_\eta \theta$ on $\mathcal{Z} \times [-\eta',\eta']$. Thus
\begin{equation}
\label{eqprf2}
\frac{|\theta(z,\eta) - \theta(z_0, \eta)|}{\eta} \leq \frac{1}{\eta}\int_0^{\eta} K|z-z_0|\,\mathrm{d}\sigma \leq  K |z-z_0|\,.
\end{equation}
We thus have from (\ref{eqprf0}),~(\ref{eqprf1})~and~(\ref{eqprf2}):
\[
|\hatt{\theta}(z,\eta) - \hatt{\theta}(z_0, 0)| \leq K |z-z_0| + |R(z_0, \eta)|\,.
\]
The first term on the RHS can be made small by choosing $|z-z_0|$ sufficiently small; the second term can be made small by choosing $\eta$ sufficiently small. This completes the proof that $\hatt{\theta}$ is continuous on $\mathcal{Z} \times \{0\}$ and hence on its entire domain. As $\hatt{\theta}$ is continuous and $D_\eta \theta(z, 0)$ is continuous (as $\theta$ is $C^1$), $R$ is continuous as the difference of two continuous functions.
\end{proof}

\section{Proof of Theorem~\ref{mainthm}}

The proof is presented with some surrounding explanation and, for readability, several subclaims are separated from the main proof into ``subproofs''. This is to allow the reader to follow the main argument without necessarily digressing into the details of each technical claim. We break the proof into numbered parts which can be referred back to.

\begin{enumerate}[align=left,leftmargin=*,itemsep=2ex,parsep=1ex]

\item 
\label{pointbasic}
{\bf The basic set-up.} We suppose that the original CRN $\mathcal{R}$ described by (\ref{eqbasic0}), namely, $\dot x = \Gamma v(x)$, admits a nondegenerate (resp., linearly stable) positive periodic orbit $\mathcal{O}$. Recall that $\Gamma$ is an $n \times r_0$ \textcolor{black}{matrix}, and that $v\colon \mathbb{R}^n_{\gg 0} \to \mathbb{R}^{r_0}$ is assumed to be $C^2$. Define $S_{\mathcal{O}}$ to be the coset of $\mathrm{im}\,\Gamma$ which contains $\mathcal{O}$, namely $S_{\mathcal{O}} = x_0+\mathrm{im}\,\Gamma$ for some $x_0 \in \mathcal{O}$. Let $\mathcal{Z}^o$ be some connected subset of $S_{\mathcal{O}}$ containing $\mathcal{O}$, relatively open w.r.t. $S_{\mathcal{O}}$, and whose closure, $\mathcal{Z}$, is compact and lies in $\mathbb{R}^n_{\gg 0}$. We do not need to introduce local coordinates on $S_{\mathcal{O}}$ explicitly, although this can be done (see Remark~\ref{remlocalcoords}). \textcolor{black}{The geometry is illustrated schematically in Figure~\ref{fig0}}.

\begin{figure}[h]
\begin{center}
\begin{tikzpicture}[domain=0:4,scale=0.4]

\draw[-, line width=0.04cm] (0,0) -- (12,0);
\draw[-, line width=0.04cm] (0,0) -- (0,11.5);
\draw[-, line width=0.04cm] (0,0) -- (9.5,4.75);

\fill[color=red!40,fill opacity=1] (0,9) -- (9,4.5) -- (8.5,0) -- cycle;

\begin{scope}[xshift=2cm,yshift=1cm,scale=0.9]
\draw[->, line width=0.04cm] (5.5, 2) .. controls (6, 2) and (6.3,3) .. (6,3.5);
\draw[->, line width=0.04cm] (6, 3.5) .. controls (5.7, 4) and (4,5) .. (3.5,4.5);
\draw[->, line width=0.04cm] (3.5, 4.5) .. controls (3, 4) and (5,2) .. (5.5,2);

\end{scope}

\begin{scope}[xshift=0cm,yshift=-0.2cm,scale=1.3]
\fill[color=black, fill opacity=0.2, -] (5.5, 1.8) .. controls (6.5, 1.6) and (6.3,3) .. (6,3.5) .. controls (5.7, 4) and (4,5.3) .. (3.2,5) .. controls (2.4, 4.7) and (4.5,2) .. (5.5,1.8);

\end{scope}

\node[rotate=10,xslant=-0.8] at (6,4.2) {$\mathcal{O}$};

\node[rotate=20,xslant=-0.8] at (2.6,7.2) {$S_\mathcal{O}$};

\node[rotate=30,xslant=-0.8] at (4.7,5.7) {$\mathcal{Z}^o$};
\end{tikzpicture}
\end{center}
\caption{\label{fig0} \textcolor{black}{A schematic showing the periodic orbit $\mathcal{O}$, the positive part of $S_\mathcal{O}$, the coset of $\mathrm{im}\,\Gamma$ on which it lies, and $\mathcal{Z}^o$, a bounded, relatively open subset of $S_\mathcal{O}$ containing $\mathcal{O}$ and with positive closure. The situation illustrated is the case $n=3$, $r=2$, namely the state space has dimension $3$ and the stoichiometry classes have dimension $2$.}}
\end{figure}

Recall the definitions of the matrices $a$, $a'$, $b$, $b'$, $\alpha$ and $\beta$. Given the assumption that the $(m+k) \times m$ matrix $\beta$ has rank $m$, we can assume without loss of generality (i.e., by reordering the added species $Y$ if necessary) that $\beta = \left(\begin{array}{c}\hatt{\beta}\\\doublehat{\beta}\end{array}\right)$, where $\hatt{\beta}$ is a nonsingular $m \times m$ matrix, and $\doublehat{\beta}$ is a $k \times m$ matrix. If $k=0$, then $\doublehat{\beta}$ is empty. $\hatt{y} \in \mathbb{R}^m$, $\doublehat{y} \in \mathbb{R}^k$, $\hatt{b} \in \mathbb{R}^{m\times m}$, $\hatt{b}' \in \mathbb{R}^{m\times m}$, $\doublehat{b} \in \mathbb{R}^{k\times m}$ and $\doublehat{b}' \in \mathbb{R}^{k\times m}$ are defined in the natural way ($\doublehat{y}$, $\doublehat{b}$, and $\doublehat{b}'$ are empty if $k=0$). 

The stoichiometric matrix of $\mathcal{R}'$ is
\[
\Gamma':=\left(\begin{array}{cc}\Gamma&\alpha\\0&\beta\end{array}\right)\,.
\]
Let $S^*$ denote the coset of $\mathrm{im}\,\Gamma'$ which includes the point $(x, \hatt{y}, \doublehat{y}) = (x_0, 0,\mathbf{1})$. (In the case $k=0$, this just means that $(x, \hatt{y}) = (x_0, 0)$.) Let $S^*_+$ denote the positive part of $S^*$, namely $S^*_+:= S^* \cap (\mathbb{R}^n_{\gg 0} \times \mathbb{R}^{m+k}_{\gg 0})$. \textcolor{black}{$S^*_+$ is the positive stoichiometry class of interest to us:} our goal is to show that we can choose rates for the added reactions from the class of mass action kinetics such that $\mathcal{R}'$ admits a periodic orbit on $S^*_+$ which is nondegenerate (resp., linearly stable) relative to $S^*$.

\item 
\label{pointcoords}
{\bf A coordinate transformation to simplify the geometry.} With the ultimate aim of setting up a singular perturbation problem \textcolor{black}{amenable to the techniques of geometric singular perturbation theory}, we now carry out a number of steps. The first is to define a new variable $z = x-\alpha\hatt{\beta}^{-1}\hatt{y}$. More precisely, we define the linear bijection $\phi_0\colon \mathbb{R}^{n} \times \mathbb{R}^{m+k} \to \mathbb{R}^{n} \times \mathbb{R}^{m+k}$ by
\[
\phi_0\left(\begin{array}{c}x\\y\end{array}\right) = \left(\begin{array}{c}x-\alpha\hatt{\beta}^{-1}\hatt{y}\\y\end{array}\right), \quad \mbox{so that} \quad \phi_0^{-1}\left(\begin{array}{c}z\\y\end{array}\right) = \left(\begin{array}{c}z+\alpha\hatt{\beta}^{-1}\hatt{y}\\y\end{array}\right)\,.
\]
We refer to the domain of $\phi_0$ as $(x,y)$-space and its codomain as $(z,y)$-space. It is easily shown that $\phi_0$ takes cosets of 
\[
\textcolor{black}{\mathrm{im}}\,\left(\begin{array}{cc}\Gamma&\alpha\\0&\beta\end{array}\right) \quad \mbox{in $(x,y)$-space to cosets of } \textcolor{black}{\mathrm{im}}\,\left(\begin{array}{cc}\Gamma&0\\0&\beta\end{array}\right)\quad \mbox{in $(z,y)$-space.}
\]
(See Subproof~\ref{phi0cosets}.) \textcolor{black}{This transformation has permitted us to write the initial stoichiometric matrix in block diagonal form at the cost of somewhat complicating the subsequent rate functions. We will see shortly that this also allows us to formulate a singular perturbation problem.} In $(z,y)$-space, (\ref{perteq}) becomes 
\begin{equation}
\label{eqzyfull}
\begin{array}{rcl}\dot z & = & \Gamma v(z+\alpha\hatt{\beta}^{-1}\hatt{y})\,,\\
\dot y & = & \beta q(z+\alpha\hatt{\beta}^{-1}\hatt{y},y)\,,
\end{array}
\end{equation}
(\ref{eqzyfull}) defines a local flow on $\phi_0(\mathbb{R}^n_{\gg 0} \times \mathbb{R}^{m+k}_{\gg 0}) = \{(z,y) \in \mathbb{R}^n \times \mathbb{R}^{m+k}\colon z +\alpha\hatt{\beta}^{-1}\hat y \gg 0,\,\, y \gg 0\}$.

\item 
\label{pointreduce}
{\bf Removing the variables $\doublehat{y}$.} We now carry out a further transformation to eliminate the variable $\doublehat{y}$ from (\ref{eqzyfull}). \textcolor{black}{The goal is to remove complications arising if the new reactions introduce new conservation laws. Conserved quantities involving only the new species can immediately be fixed at arbitrary positive values, eliminating some added species and simplifying later arguments. A similar strategy was adopted in the proof of Theorem~5 of \cite{banajipanteaMPNE}.} Define $\delta := -(\doublehat{\beta}\hatt{\beta}^{-1})^{\mathrm{t}}$. Then $\delta^{\mathrm{t}}\hatt{y}+\doublehat{y}$ is constant along trajectories of (\ref{eqzyfull}). (See Subproof~\ref{yhathat}.) We fix the value of $\delta^{\mathrm{t}}\hatt{y}+\doublehat{y}$ as $\mathbf{1}$ (this choice is arbitrary: any fixed vector in $\mathbb{R}^k_{\gg 0}$ would do) and define the hyperplane
\[
\mathcal{H}:=\{(z,y) \in \mathbb{R}^n \times \mathbb{R}^{m+k}\colon \delta^{\mathrm{t}}\hatt{y}+\doublehat{y}=\mathbf{1}\}\,.
\]
Note that according to our conventions, $\mathcal{H}:=\mathbb{R}^n \times \mathbb{R}^{m}$ if $k=0$. \textcolor{black}{Trajectories of (\ref{eqzyfull}) beginning on $\mathcal{H}$ remain on $\mathcal{H}$, as $\mathcal{H}$ is} a union of cosets of $\mathrm{im}\,\Gamma \times \mathrm{im}\beta$ in $\mathbb{R}^n \times \mathbb{R}^{m+k}$. Define the affine bijection $\phi_1\colon \mathcal{H} \to \mathbb{R}^n \times \mathbb{R}^m$ by $\phi_1(z,(\hatt{y},\mathbf{1}-\delta^{\mathrm{t}}\hatt{y}))=(z, \hatt{y})$. $\phi_1$ is the projection $(z, \hatt{y}, \doublehat{y}) \mapsto (z, \hatt{y})$ restricted to $\mathcal{H}$ and is just the identity on $\mathbb{R}^n \times \mathbb{R}^m$ if $k=0$. Define $S:=S_{\mathcal{O}} \times \mathbb{R}^m$, and observe that 
\begin{equation}
\label{eqsodash}
\phi_0^{-1}\circ \phi_1^{-1}(S) = S^*\,.
\end{equation}
(See Subproof~\ref{phi0phi1}.) \textcolor{black}{Also define $\mathcal{H}' = \phi_0^{-1}(\mathcal{H})$, namely $\mathcal{H}'$ is the subset of $(x,y)$-space defined by the equation $\delta^{\mathrm{t}}\hatt{y}+\doublehat{y}=\mathbf{1}$.} The action of the transformations $\phi_0$ and $\phi_1$ is summarised in Figure~\ref{fig3}.

\begin{figure}[h]
\begin{center}
\begin{tikzpicture}
\draw[thick] (0,2) .. controls (0,3.8) and (3.6,4) .. (4,2) -- (4,2) .. controls (4.4,0) and (0,0) .. (0,2);
\node at (2,3) {$(x,y)$-space};
\fill[black!20, opacity=0.5] (0.6,1.5) .. controls (0.6,2.7) and (1.8,2.7) .. (2,2.7) -- (2,2.7) .. controls (2.5,2.7) and (3.3,2.7) .. (3.3,1.5) -- (3.3,1.5) .. controls (3.3,1) and (2.5,1) .. (2,1) -- (2,1) .. controls (1.5,1) and (0.6,1) .. (0.6,1.5);
\draw[thick] (0.6,1.5) .. controls (0.6,2.7) and (1.8,2.7) .. (2,2.7) -- (2,2.7) .. controls (2.5,2.7) and (3.3,2.7) .. (3.3,1.5) -- (3.3,1.5) .. controls (3.3,1) and (2.5,1) .. (2,1) -- (2,1) .. controls (1.5,1) and (0.6,1) .. (0.6,1.5);
\node at (2,2.2) {$\delta^{\mathrm{t}}\hatt{y}+\doublehat{y}=\mathbf{1}$};
\node at (1,1.5) {$\mathcal{H}'$};

\fill[black!30, opacity=0.5] (2.2,1.5) ellipse (0.6cm and 0.3cm);
\draw (2.2,1.5) ellipse (0.6cm and 0.3cm);
\node at (2.2,1.5) {$S^*$};

\begin{scope}[xshift=5cm]
\draw[thick] (0,2) .. controls (0,3.8) and (3.6,4) .. (4,2) -- (4,2) .. controls (4.4,0) and (0,0) .. (0,2);
\node at (2,3) {$(z,y)$-space};
\fill[black!20, opacity=0.5] (0.6,1.5) .. controls (0.6,2.7) and (1.8,2.7) .. (2,2.7) -- (2,2.7) .. controls (2.5,2.7) and (3.3,2.7) .. (3.3,1.5) -- (3.3,1.5) .. controls (3.3,1) and (2.5,1) .. (2,1) -- (2,1) .. controls (1.5,1) and (0.6,1) .. (0.6,1.5);
\draw[thick] (0.6,1.5) .. controls (0.6,2.7) and (1.8,2.7) .. (2,2.7) -- (2,2.7) .. controls (2.5,2.7) and (3.3,2.7) .. (3.3,1.5) -- (3.3,1.5) .. controls (3.3,1) and (2.5,1) .. (2,1) -- (2,1) .. controls (1.5,1) and (0.6,1) .. (0.6,1.5);
\node at (2,2.2) {$\delta^{\mathrm{t}}\hatt{y}+\doublehat{y}=\mathbf{1}$};
\node at (2,1.5) {$\mathcal{H}$};

\end{scope}

\begin{scope}[xshift=10cm]
\fill[black!20, opacity=0.5] (0.6,1.5) .. controls (0.6,2.7) and (1.8,2.7) .. (2,2.7) -- (2,2.7) .. controls (2.5,2.7) and (3.3,2.7) .. (3.3,1.5) -- (3.3,1.5) .. controls (3.3,1) and (2.5,1) .. (2,1) -- (2,1) .. controls (1.5,1) and (0.6,1) .. (0.6,1.5);
\draw[thick] (0.6,1.5) .. controls (0.6,2.7) and (1.8,2.7) .. (2,2.7) -- (2,2.7) .. controls (2.5,2.7) and (3.3,2.7) .. (3.3,1.5) -- (3.3,1.5) .. controls (3.3,1) and (2.5,1) .. (2,1) -- (2,1) .. controls (1.5,1) and (0.6,1) .. (0.6,1.5);
\node at (2,2.2) {$(z,\hatt{y})$-space};

\fill[black!30, opacity=0.5] (2.2,1.5) ellipse (0.6cm and 0.3cm);
\draw (2.2,1.5) ellipse (0.6cm and 0.3cm);
\node at (2.2,1.5) {$S$};

\end{scope}

\draw[thick,-latex] (3.2,3.2) .. controls (3.8,3.7) and (5.2,3.7) .. (5.8,3.2);
\node at (4.5,3.8) {$\phi_0$};

\draw[thick,-latex] (3.2,2.2) .. controls (3.8,2.7) and (5.2,2.7) .. (5.8,2.2);
\node at (4.5,2.8) {$\phi_0$};

\draw[thick,-latex] (8.2,2.2) .. controls (8.8,2.7) and (10.2,2.7) .. (10.8,2.2);
\node at (9.5,2.8) {$\phi_1$};

\draw[thick,latex-] (11.8,1.2) .. controls (9.8,-0.2) and (4.5,-0.2) .. (2.5,1.2);
\node at (7.1,-0.15) {$\phi_1 \circ \phi_0$};

\end{tikzpicture}
\end{center}
\caption{\label{fig3} \textcolor{black}{The bijection $\phi_0$ takes $(x,y)$-space to $(z,y)$-space and $\mathcal{H}'$ to $\mathcal{H}$. The bijection $\phi_1$ takes $\mathcal{H}$ to $(z, \hatt{y})$-space.} $\phi_1\circ \phi_0$ is defined on $\textcolor{black}{\mathcal{H}'}$, and is an affine bijection between $\mathcal{H}'$ and $(z, \hatt{y})$-space. Its restriction to $S^*$ is an affine bijection between $S^*$ and $S$.}
\end{figure}

The claim that $\phi_0^{-1}\circ \phi_1^{-1}(S) = S^*$ is consistent with the following notational convention which we now adopt: given any set $X$ in $(z,\hatt{y})$-space, $X^{*}$ refers to the corresponding set in $(x,y)$-space, namely $X^{*}=\phi_0^{-1}\circ \phi_1^{-1}(X)$. Notice that if $X \subseteq S$, then $X^* \subseteq S^*$. In the other direction, given a function $F$ on $(x,y)$-space, $F_{*}$ refers to the corresponding function on $(z,\hatt{y})$-space, namely $F_{*}=F\circ\phi_0^{-1}\circ \phi_1^{-1}$.

\item {\bf The system in $(z,\hatt{y})$-space.} 
\label{zysys}
Define \textcolor{black}{the following subsets of $(z,\hatt{y})$-space.}
\[
W_+:=\{(z,\hatt{y}) \in \mathbb{R}^n \times \mathbb{R}^m \colon \phi_0^{-1}\circ \phi_1^{-1}(z, \hatt{y}) \in \mathbb{R}^n_{\gg 0} \times \mathbb{R}^{m+k}_{\gg 0}\},
\] 
or explicitly $W_+=\{(z,\hatt{y}) \in \mathbb{R}^n \times \mathbb{R}^m \colon z+\alpha\hatt{\beta}^{-1}\hatt{y} \gg 0,\,\, \hatt{y} \gg 0,\,\, \mathbf{1}-\delta^{\mathrm{t}}\hatt{y} \gg 0\}$. Similarly, define
\[
S_+:=\{(z,\hatt{y}) \in S \colon \phi_0^{-1}\circ \phi_1^{-1}(z, \hatt{y}) \in \mathbb{R}^n_{\gg 0} \times \mathbb{R}^{m+k}_{\gg 0}\},
\]
namely, $S_+ := \{(z,\hatt{y}) \in S \colon z+\alpha\hatt{\beta}^{-1}\hatt{y} \gg 0,\,\, \hatt{y} \gg 0,\,\, \mathbf{1}-\delta^{\mathrm{t}}\hatt{y} \gg 0\}$. Note that $S_+^*= \phi_0^{-1}\circ \phi_1^{-1}(S_+)$ as expected from our notational conventions.

\textcolor{black}{$W_+$ is the image in $(z,\hatt{y})$-space of the positive part of $\mathcal{H}'$, namely, $W_+$ is the image under $\phi_1 \circ \phi_0$ of a union of positive stoichiometry classes; and $S_+$ is the image under $\phi_1 \circ \phi_0$ in $(z,\hatt{y})$-space of the positive part of $S^*$, the positive stoichiometry class of interest to us.}

Define $q_{*} = q\circ \phi_0^{-1}\circ\phi_1^{-1}$, namely, $q_{*}(z,\hatt{y}) = q(z+\alpha\hatt{\beta}^{-1}\hatt{y},(\hatt{y}, \mathbf{1}-\delta^{\mathrm{t}}\hatt{y}))$. As $q$ is defined on $\mathbb{R}^n_{\gg 0} \times \mathbb{R}^{m+k}_{\gg 0}$, $q_{*}$ is defined on $W_+$. The restriction of (\ref{eqzyfull}) to $\mathcal{H}$, followed by projection by $\phi_1$, gives the following system on $W_+$:
\begin{equation}
\label{eqzyhat1}
\begin{array}{rcl}\dot z & = & \Gamma v(z+\alpha\hatt{\beta}^{-1}\hatt{y})\,,\\
\dot {\hatt{y}} & = & \hatt{\beta} q_{*}(z,\hatt{y})\,.
\end{array}
\end{equation}
We want to restrict attention to a \textcolor{black}{region of $S_+$ on which the values of $\hatt{y}$ are small and positive. To this end,} we observe that there exists $y_{pos}>0$ such that
\[
\mathcal{Z}_+:=\{(z,\hatt{y})\colon z \in \mathcal{Z},\,\, \hatt{y} \gg 0, |\hatt{y}| \leq y_{pos}\}
\]
lies in $S_+$, namely,
\begin{equation}
\label{eqmapstopos}
\mathcal{Z}_+^*:=\phi_0^{-1}\circ \phi_1^{-1}(\mathcal{Z}_+) \subseteq \mathbb{R}^n_{\gg 0} \times \mathbb{R}^{m+k}_{\gg 0}\,.
\end{equation}
(See Subproof~\ref{zplus}.) Since $\mathcal{Z}_+ \subseteq S_+$, it follows that $\mathcal{Z}_+^* \subseteq S^*_+$. Thus (\ref{eqzyhat1}) defines a local flow on $\mathcal{Z}_+^o$, the relative interior of $\mathcal{Z}_+$ in $S$, and applying $\phi_0^{-1}\circ \phi_1^{-1}$ we get a corresponding local flow on a relatively open subset of $S^*_+$. \textcolor{black}{If $\hatt{y} \in \mathbb{R}$ for example, then $\mathcal{Z}_+ = \mathcal{Z} \times (0, y_{pos}]$, and $\mathcal{Z}_+^o = \mathcal{Z}^o \times (0, y_{pos})$. The eventual goal of our constructions is to prove the existence of a periodic orbit of (\ref{eqzyhat1}) in $\mathcal{Z}_+^o$ for appropriate choices of $q_{*}$.}

\item 
\label{pointrates}
{\bf Choosing the rates for the new reactions.} The rate function $q$ of the new reactions has so far been left undetermined. $q$ will be chosen from some class of rate functions depending on two \textcolor{black}{real} parameters $\epsilon$ and $\eta$, and will take the form
\[
\frac{1}{\epsilon}f(x,y,\eta)\,.
\]
$f$ will be chosen shortly, but for the moment we assume that $f$ is defined and $C^2$ on $\mathbb{R}^n_{\gg 0} \times \mathbb{R}^{m+k}_{\gg 0} \times \mathbb{R}_{> 0}$. The reason for introducing the two parameters is that, roughly speaking, we need to be able to make both rates of each reversible reaction arbitrarily large (via $\epsilon$) while independently controlling the ratio of forward and backward rates (via $\eta$). \textcolor{black}{The intuition is that (i) we want to be able to ensure that the added reactions are close to equilibrium via $\epsilon$ (this will be our singular perturbation parameter); and (ii) we want to ensure, via $\eta$, that the values of $\hatt{y}$ are small when the added reactions are at equilibrium.}

In arguments below we sometimes extend $\phi_0$ and $\phi_1$ to $(x,y,\eta)$-space and $(z,y,\eta)$-space respectively, namely we let $\phi_0$ and $\phi_1$ refer to $\phi_0 \times \mathrm{id}$ and $\phi_1 \times \mathrm{id}$. This should cause no confusion. 

Define $f_{*}:=f\circ \phi_0^{-1}\circ\phi_1^{-1}$, or more explicitly $f_{*}(z,\hatt{y}, \eta) = f(z+\alpha\hatt{\beta}^{-1}\hatt{y},(\hatt{y}, \mathbf{1}-\delta^{\mathrm{t}}\hatt{y}),\eta)$. By construction $f_{*}$ is defined and $C^2$ on $W_+\times \mathbb{R}_{>0}$ which includes $\mathcal{Z}_+ \times \mathbb{R}_{> 0}$. (\ref{eqzyhat1}) now takes the form of a typical singular perturbation problem:
\begin{equation}
\label{eqzyhat1a}
\begin{array}{rcl}\dot z & = & \Gamma v(z+\alpha\hatt{\beta}^{-1}\hatt{y})\,,\\
\epsilon\dot {\hatt{y}} & = & \hatt{\beta} f_{*}(z,\hatt{y},\eta)\,.
\end{array}
\end{equation}
When we study (\ref{eqzyhat1a}) we restrict attention to $\mathcal{Z}_+^o$.

We have left fixing of the kinetics of the added reactions to this late stage in order to facilitate generalisation of the results. But now we assume mass action kinetics for the added reactions (\ref{addedreacs}), and set
\begin{equation}
\label{eqf}
 f(x,y, \eta) := \bm{\eta}^{-\hatt{b}^{\mathrm{t}}}\circ x^{a^{\mathrm{t}}}\circ y^{b^{\mathrm{t}}} - \bm{\eta}^{-{\hatt{b}'}^{\mathrm{t}}}\circ x^{{a'}^{\mathrm{t}}}\circ y^{{b'}^{\mathrm{t}}}\,.
\end{equation}

\item 
\label{pointposzero}
{\bf The positive zero-set of $f$.} We will be interested in the positive zero-set of $f$, \textcolor{black}{corresponding to the added reactions being at equilibrium}. Let $\gamma := -(\alpha\,\hatt{\beta}^{-1})^{\mathrm{t}}$ and recall that we defined $\delta = -(\doublehat{\beta}\,\hatt{\beta}^{-1})^{\mathrm{t}}$. With some manipulation we find that for each fixed $\eta>0$, solutions to $f(x,y, \eta)=0$ on $\mathbb{R}^n_{\gg 0} \times \mathbb{R}^{m+k}_{\gg 0}$ are precisely solutions to 
\[
g(x,y,\eta):=\hatt{y} - \eta x^\gamma\circ\doublehat{y}^\delta=0
\]
on $\mathbb{R}^n_{\gg 0} \times \mathbb{R}^{m+k}_{\gg 0}$ (see Subproof~\ref{ftog}.) Note that, for any $\gamma$ and $\delta$, $g$ is defined on \textcolor{black}{a larger domain than $f$, namely,} $\mathbb{R}^n_{\gg 0} \times (\mathbb{R}^m\times \mathbb{R}^k_{\gg 0}) \times \mathbb{R}$ (this is to be interpreted as $\mathbb{R}^n_{\gg 0} \times \mathbb{R}^m \times \mathbb{R}$ if $k=0$). Define $g_{*} := g\circ\phi_0^{-1}\circ\phi_1^{-1}$, namely,
\[
g_{*}(z,\hatt{y},\eta)=\hatt{y} - \eta (z+\alpha\hatt{\beta}^{-1}\hatt{y})^\gamma\circ (\mathbf{1}-\delta^{\mathrm{t}}\hatt{y})^\delta\,.
\]
$g_{*}$ is defined and $C^2$ provided $z+\alpha\hatt{\beta}^{-1}\hatt{y} \gg 0$ and $\mathbf{1}-\delta^{\mathrm{t}}\hatt{y} \gg 0$, and so the domain of $g_{*}$ includes an open neighbourhood of $\mathcal{Z} \times B_{y_{pos}} \times \mathbb{R}$ (recall Part~\ref{zysys}). It is important in constructions below involving the IFT that the domain of $g_{*}$ includes $\mathcal{Z} \times \{0\} \times \{0\}$.

\item 
\label{pointfirstbound}
{\bf An upper bound on $\eta$.} We need to put a number of upper bounds on $\eta$. The first of these ensures (via the IFT) that a portion of the zero set of $g_{*}$, and hence $f_{*}$, is the graph of a function. Let $\mathcal{E}$ be the zero set of $g_{*}$ on its domain, namely
\[
\mathcal{E}:=\{(z, \hatt{y}, \eta)\in \mathbb{R}^n \times \mathbb{R}^m \times \mathbb{R}\colon z+\alpha\hatt{\beta}^{-1}\hatt{y} \gg 0,\,\, \mathbf{1}-\delta^{\mathrm{t}}\hatt{y} \gg 0,\,\, g_{*}(z,\hatt{y},\eta)=0\}\,.
\]
We now claim that there exists $y_{max}\in(0, y_{pos}]$ and $\eta_1>0$ such that the zero-set of $g_{*}$ in $\mathcal{Z} \times B_{y_{max}} \times (-\eta_1, \eta_1)$ is the graph of a $C^2$ function $\theta\colon \mathcal{Z} \times (-\eta_1,\eta_1) \to \mathbb{R}^m$, namely, 
\[
\mathcal{E} \cap (\mathcal{Z} \times B_{y_{max}} \times (-\eta_1, \eta_1)) = \{(z,\hatt{y}, \eta) \colon z \in \mathcal{Z},\,\, \eta \in (-\eta_1,\eta_1),\,\,\hatt{y} = \theta(z, \eta)\}\,.
\]
\textcolor{black}{This is shown by applying the extended IFT in Lemma~\ref{lemIFT1}} -- (see Subproof~\ref{theta}). It is also clear that $\theta$ must satisfy $\theta(z,0)=0$. For each fixed $\eta \in (-\eta_1, \eta_1)$ we then have 
\[
\mathcal{E}_\eta:=\{(z, \hatt{y}) \in \mathcal{Z} \times B_{y_{max}} \colon g_{*}(z,\hatt{y},\eta)=0\} = \{(z,\hatt{y}) \colon z \in \mathcal{Z},\,\,\hatt{y} = \theta(z, \eta)\}\,\,.
\]
The geometry of the situation is illustrated schematically in Figure~\ref{fig4}.

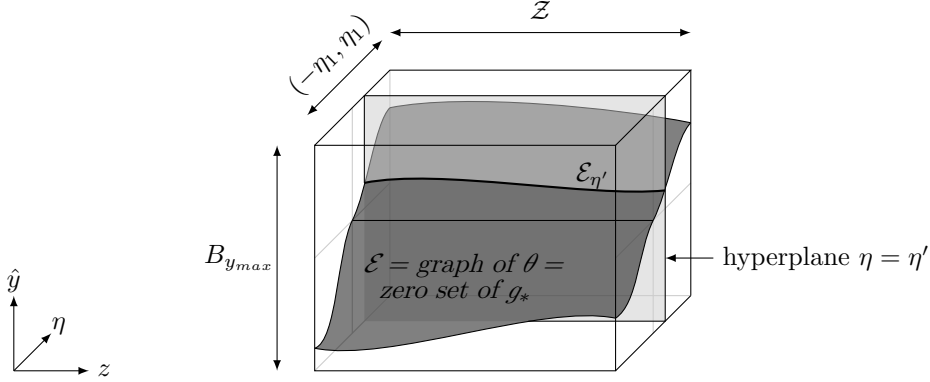
\begin{figure}[h]
\begin{center}
\begin{tikzpicture}
\draw[color=black!20] (0,0)--(1,1) -- (5,1);
\draw[color=black!20] (1,1)--(1,4);
\draw[color=black!20] (0.5,0.5)--(0.5,3.5);
\draw[color=black!20] (4.5,0.5)--(4.5,3.5);
\draw[color=black!20] (0,1.5)--(1,2.5);
\draw[color=black!20] (4,1.5)--(5,2.5);

\fill[black!20, opacity=0.5] (4.66,2.4) .. controls (3.65,2.3) and (1.65, 2.7) .. (0.66, 2.5) -- (0.66,2.5)--(0.66,0.66) -- (4.66,0.66) -- (4.66,2.4);

\draw (0.5,2) -- (4.5,2);

\fill[black,opacity=0.5] (0.5,2) .. controls (0.7,2.4) and (0.8,3.3) .. (1,3.5) -- (1,3.5) .. controls (2,3.7) and (4,3.5) .. (5,3.3) -- (5,3.3) .. controls (4.8,3.1) and (4.7,2.4) .. (4.5,2) -- (4.5,2) .. controls (4.3,1.6) and (4.3,0.9) .. (4,0.7) -- (4,0.7) .. controls (3,0.9) and (1,0.1) .. (0,0.3) -- (0,0.3) .. controls (0.3,0.4) and (0.3,1.6) .. (0.5,2);
\draw (0.5,2) .. controls (0.7,2.4) and (0.8,3.3) .. (1,3.5);
\draw (1,3.5) .. controls (2,3.7) and (4,3.5) .. (5,3.3);
\draw (5,3.3) .. controls (4.8,3.1) and (4.7,2.4) .. (4.5,2);
\draw (4.5,2) .. controls (4.3,1.6) and (4.3,0.9) .. (4,0.7);
\draw (4,0.7) .. controls (3,0.9) and (1,0.1) .. (0,0.3);
\draw (0,0.3) .. controls (0.3,0.4) and (0.3,1.6) .. (0.5,2);

\fill[black!20, opacity=0.5] (0.66,2.5)--(0.66,3.66) -- (4.66,3.66) -- (4.66,2.4) -- (4.66,2.4) .. controls (3.65,2.3) and (1.65, 2.7) .. (0.66, 2.5);
\draw (0.66,2.5)--(0.66,3.66) -- (4.66,3.66) -- (4.66,2.4) -- (4.66,0.66) -- (2.84,0.66);
\draw[thick] (0.66,2.5) .. controls (1.65,2.7) and (3.65,2.3) .. (4.66,2.4);

\draw (0,0)--(4,0) -- (4,3) -- (0,3) -- cycle;
\draw (0,3)--(1,4) -- (5,4) -- (4,3);
\draw (5,1) -- (4,0);
\draw (5,1)--(5,4);
\draw[latex-latex] (1,4.5) -- (5,4.5);
\node at (3,4.8) {$\mathcal{Z}$};
\draw[latex-latex] (-0.1,3.4) -- (0.9,4.4);
\node[rotate=45] at (0.2,4.2) {$(-\eta_1,\eta_1)$};
\draw[latex-latex] (-0.5,0) -- (-0.5,3);
\node at (-1,1.5) {$B_{y_{max}}$};
\begin{scope}[cm={1,0,0.3,1,(2,1.4)}]
\node[transform shape] at (0,0) {$\mathcal{E} = $ graph of $\theta = $};
\node[transform shape] at (0,-0.35) {zero set of $g_*$};
\end{scope}

\draw[-latex] (5.3,1.5) -- (4.66,1.5);
\node at (6.8,1.5) {hyperplane $\eta=\eta'$};

\node at (3.7,2.6) {$\mathcal{E}_{\eta'}$};

\draw[-latex] (-4,0)--(-3,0);
\draw[-latex] (-4,0)--(-4,1);
\draw[-latex] (-4,0)--(-3.5,0.5);
\node at (-2.8,0) {$z$};
\node at (-4,1.2) {$\hatt{y}$};
\node at (-3.4,0.6) {$\eta$};

\end{tikzpicture}
\end{center}
\caption{\label{fig4} The graph of $\theta$ coincides with the zero set of $g_*$ inside $\mathcal{Z} \times B_{y_{max}} \times (-\eta_1,\eta_1)$. By choosing $\eta_1$ to be sufficiently small we can ensure that $y_{max}$ is as small as we like. For each $\eta' \in (-\eta_1, \eta_1)$, $\mathcal{E}_{\eta'}$ is the intersection of this graph with the hyperplane $\eta=\eta'$. \textcolor{black}{Note that this figure is very much a schematic: in practice $\mathcal{Z}$, and hence each set $\mathcal{E}_{\eta}$, must be at least two dimensional.}}
\end{figure}

\item 
\label{pointsecondbound}
{\bf A second upper bound on $\eta$.} This is needed to ensure that the function $\theta$ just constructed is strictly positive for positive $\eta$. In other words, we claim that there exists $\eta_2\in (0,\eta_1]$, such that $(z,\eta) \in \mathcal{Z} \times (0,\eta_2)$ implies $\theta(z,\eta) \gg 0$ (see Subproof~\ref{ypos}.) Recall, additionally, (Part~\ref{pointfirstbound} above) that $(z,\eta) \in \mathcal{Z} \times [0,\eta_1)$ implies that $|\theta(z,\eta)|<y_{max}$. Thus $(z,\eta) \in \mathcal{Z} \times (0,\eta_2)$ implies that $(z, \theta(z,\eta)) \in \mathcal{Z}_+$. As a result, provided $\eta \in (0,\eta_2)$, $\mathcal{E}_\eta \subseteq \mathcal{Z}_+\subseteq S_+$ (see Part~\ref{zysys}), and consequently, defining $\mathcal{E}_\eta^{*}:=\phi_0^{-1}\circ \phi_1^{-1}(\mathcal{E}_\eta)$, \textcolor{black}{we have,}
\begin{equation}
\label{eqmapstopos1}
\mathcal{E}_\eta^{*} \subseteq \mathcal{Z}_+^* \subseteq S_+^*\,.
\end{equation}

\item 
\label{pointthirdbound}
{\bf A third upper bound on $\eta$.} This is needed to ensure that the differential algebraic system obtained in a singular limit has a nondegenerate (resp., linearly stable) periodic orbit. Fix any $\eta \in (0,\eta_2)$ and consider the following system
\begin{equation}
\label{eqredZ}
\dot z = \Gamma v(z+\alpha\hatt{\beta}^{-1}\theta(z, \eta))\,,
\end{equation}
on $\mathcal{Z}^o$. \textcolor{black}{The reader may glance ahead to Part~\ref{pointSPT} of the proof and note that} the vector field $F(z, \eta):=\Gamma v(z+\alpha\hatt{\beta}^{-1}\theta(z, \eta))$ occurring on the right of (\ref{eqredZ}) can be regarded as the so-called ``reduced vector field'' associated with (\ref{eqzyhat1a}). Note that $F$ is $C^2$ on $\mathcal{Z} \times (0,\eta_2)$ and, by assumption, $\dot z = F(z,0) = \Gamma v(z)$ has a periodic orbit $\mathcal{O}$ in $\mathcal{Z}^o$, nondegenerate relative to $S_{\mathcal{O}}$. By regular perturbation theory arguments (Lemma~\ref{lemreg}) there exists $\eta_3\in(0,\eta_2]$ such that provided $\eta\in(0,\eta_3)$ (\ref{eqredZ}) has a periodic orbit $\mathcal{O}'_\eta$  in $\mathcal{Z}^o$ close to $\mathcal{O}$, and such that the number of Floquet multipliers of $\mathcal{O}'_\eta$ relative to $S_{\mathcal{O}}$ inside and outside the unit circle is the same as that of $\mathcal{O}$. Consequently if $\mathcal{O}$ is nondegenerate (resp., linearly stable) relative to $S_{\mathcal{O}}$, then so is $\mathcal{O}'_\eta$.

\item 
\label{pointfourthbound}
{\bf A fourth upper bound on $\eta$.} We need one more upper bound on $\eta$, connected with ensuring that normal hyperbolicity conditions needed to apply results in \cite{Fenichel79} hold. For each $(z,\eta) \in \mathcal{Z} \times (0, \eta_3)$ define $W(z,\eta) := \hatt{\beta} D_{\hatt{y}}f_{*}(z, \theta(z,\eta), \eta)$, namely $W(z,\eta)$ is the Jacobian matrix of $\hatt{\beta}f_{*}(z, \hatt{y}, \eta)$ w.r.t. $\hatt{y}$, evaluated at $\mathcal{E}_\eta$. Since $f_{*}$ is defined and $C^2$ on $\mathcal{Z}_+ \times \mathbb{R}_{> 0}$, $D_{\hatt{y}}f_{*}$ is defined and $C^1$ on $\mathcal{Z}_+ \times \mathbb{R}_{> 0}$. Since, additionally, $\theta$ is $C^2$ on $\mathcal{Z} \times (0,\eta_3)$, $W$ is defined and $C^1$ on $\mathcal{Z} \times (0,\eta_3)$. We claim that there exists $\eta_4 \in (0, \eta_3]$, such that $(z,\eta) \in \mathcal{Z} \times (0,\eta_4)$ implies that $W(z,\eta)$ is Hurwitz stable, namely the eigenvalues of $W(z,\eta)$ lie in the open left half of the complex plane. The calculations are fairly lengthy and are presented in Subproof~\ref{normhyp}.

\item
\label{pointSPT}
 {\bf Singular perturbation theory: completing the argument.} We have done the preliminary work and are ready to apply perturbation theory results of Fenichel \cite{Fenichel79}. We fix some $\eta \in (0,\eta_4)$, and return to system (\ref{eqzyhat1a}) on $\mathcal{Z}_+^o$, and the equivalent ``fast time'' system obtained by rescaling time. These are
\[
\mbox{($\mathrm{A}_\epsilon$)}\,\,\,\begin{array}{rcl}\dot z & = & \Gamma v(z+\alpha\hatt{\beta}^{-1}\hatt{y})\,,\\
\epsilon\dot {\hatt{y}} & = & \hatt{\beta} f_{*}(z,\hatt{y},\eta)\,,
\end{array} \quad \mbox{and} \quad \mbox{($\mathrm{B}_\epsilon$)}\,\,\,\begin{array}{rcl}\dot z & = & \epsilon \Gamma v(z+\alpha\hatt{\beta}^{-1}\hatt{y})\,,\\
\dot {\hatt{y}} & = & \hatt{\beta} f_{*}(z,\hatt{y},\eta),
\end{array}
\]
with their respective limiting systems in the limit $\epsilon \to 0+$:
\[
\mbox{($\mathrm{A_0}$)}\,\,\,\begin{array}{rcl}\dot z & = & \Gamma v(z+\alpha\hatt{\beta}^{-1}\theta(z, \eta))\,,
\end{array} \quad \mbox{and} \quad \mbox{($\mathrm{B_0}$)}\,\,\,\begin{array}{rcl}\dot z & = & 0,\\
\dot {\hatt{y}} & = & \hatt{\beta} f_{*}(z,\hatt{y},\eta)\,.
\end{array}
\]
At this point, it may help to visit Remark~\ref{remlocalcoords}: we could, if we chose, identify $S_{\mathcal{O}}$ with $\mathbb{R}^r$ and hence $S=S_{\mathcal{O}} \times \mathbb{R}^m$ with $\mathbb{R}^{r+m}$ via a linear change of coordinates, in which case $\mathcal{Z}^o$ would be identified with an open subset of $\mathbb{R}^r$ and $\mathcal{Z}_+^o$ would be identified with an open subset of $\mathbb{R}^{r+m}$. However, instead, we choose to keep our current coordinate system and bear in mind that we are restricting attention to the locally invariant set $\mathcal{Z}_+^o$. Claims about ``nondegeneracy'' and ``linear stability'' are relative to $S_{\mathcal{O}}$ or $S$ as appropriate. 

(i) Note that $\mathcal{E}_\eta^o := \{(z,\hatt{y}) \colon z \in \mathcal{Z}^o,\,\,\hatt{y} = \theta(z, \eta)\}$ is an invariant manifold of ($\mathrm{B_0}$) consisting entirely of equilibria. Our computation showing that $W(z,\eta) = \hatt{\beta} D_{\hatt{y}}f_{*}(z, \theta(z,\eta), \eta)$ is Hurwitz (Part~\ref{pointfourthbound} above) means that the eigenvalues associated with these equilibria corresponding to directions within $S$ but transverse to $\mathcal{E}_\eta^o$ all have negative real parts. The situation is illustrated schematically in Figure~\ref{fig5}.

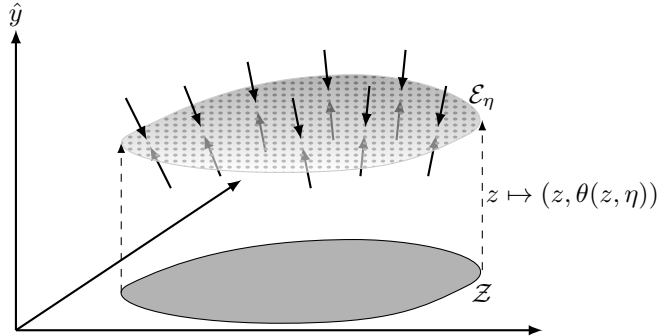
\begin{figure}[h]
\begin{center}
\begin{tikzpicture}
\draw[thick,-latex] (0,0) -- (7,0);
\draw[thick,-latex] (0,0) -- (0,4);
\node at (0,4.2) {$\hatt{y}$};
\draw[thick,-latex] (0,0) -- (3,2);

\draw[thick] (1.5,0.6) .. controls (2.1,0.9) and (3,1.2) .. (4.5,1.2) -- (4.5,1.2) .. controls (5.5,1.2) and (6.6,0.9) .. (6,0.6) -- (6,0.6) .. controls (5.4,0.3) and (5,0.1) ..  (3.5,0.1) -- (3.5,0.1) .. controls (2,0.1) and (1.1,0.4) .. (1.5,0.6);
\fill[black!30] (1.5,0.6) .. controls (2.1,0.9) and (3,1.2) .. (4.5,1.2) -- (4.5,1.2) .. controls (5.5,1.2) and (6.6,0.9) .. (6,0.6) -- (6,0.6) .. controls (5.4,0.3) and (5,0.1) ..  (3.5,0.1) -- (3.5,0.1) .. controls (2,0.1) and (1.1,0.4) .. (1.5,0.6);
\node at (6.2,0.5) {$\mathcal{Z}$};

\draw[-latex, dashed] (6.2,0.8) -- (6.2,2.8);
\draw[-latex, dashed] (1.4,0.5) -- (1.4,2.5);
\node at (7.4,1.8) {$z\mapsto (z,\theta(z,\eta))$};

\begin{scope}[yshift=2cm]
\draw[thick,latex-] (0.8+0.48*2+0.02,0.01+0.32*1.5-0.04) -- (0.8+0.48*2+0.3,0.01+0.32*1.5-0.6);
\draw[thick,latex-] (0.8+0.48*3.5+0.02,0.01+0.32*2-0.05) -- (0.8+0.48*3.5+0.24,0.01+0.32*2-0.6);
\draw[thick,latex-] (0.8+0.48*5+0.01,0.01+0.32*3-0.05) -- (0.8+0.48*5+0.12,0.01+0.32*3-0.6);
\draw[thick,latex-] (0.8+0.48*6.25+0.01,0.01+0.32*1.5-0.05) -- (0.8+0.48*6.25+0.12,0.01+0.32*1.5-0.6);
\draw[thick,latex-] (0.8+0.48*7+0.005,0.01+0.32*3.5-0.05) -- (0.8+0.48*7+0.06,0.01+0.32*3.5-0.6);
\draw[thick,latex-] (0.8+0.48*8-0.005,0.01+0.32*2-0.05) -- (0.8+0.48*8-0.06,0.01+0.32*2-0.6);
\draw[thick,latex-] (0.8+0.48*9-0.005,0.01+0.32*3.5-0.05) -- (0.8+0.48*9-0.06,0.01+0.32*3.5-0.6);
\draw[thick,latex-] (0.8+0.48*10-0.01,0.01+0.32*2-0.05) -- (0.8+0.48*10-0.12,0.01+0.32*2-0.6);

\begin{scope}
\clip  (1.5,0.6) .. controls (2.1,0.9) and (3,1.4) .. (4.5,1.4) -- (4.5,1.4) .. controls (5.5,1.4) and (6.6,0.9) .. (6,0.6) -- (6,0.6) .. controls (5.4,0.3) and (5,0.1) ..  (3.5,0.1) -- (3.5,0.1) .. controls (2,0.1) and (1.1,0.4) .. (1.5,0.6);
\foreach \p in {1,...,44}{
\foreach \q in {1,...,20}
{ 
\fill[black!60] (0.8+0.12*\p,0.01+0.08*\q) ellipse (0.03 and 0.02);
}}
\end{scope}

\shade[opacity=0.6] (1.5,0.6) .. controls (2.1,0.9) and (3,1.4) .. (4.5,1.4) -- (4.5,1.4) .. controls (5.5,1.4) and (6.6,0.9) .. (6,0.6) -- (6,0.6) .. controls (5.4,0.3) and (5,0.1) ..  (3.5,0.1) -- (3.5,0.1) .. controls (2,0.1) and (1.1,0.4) .. (1.5,0.6);
\draw[black!20] (1.5,0.6) .. controls (2.1,0.9) and (3,1.4) .. (4.5,1.4) -- (4.5,1.4) .. controls (5.5,1.4) and (6.6,0.9) .. (6,0.6) -- (6,0.6) .. controls (5.4,0.3) and (5,0.1) ..  (3.5,0.1) -- (3.5,0.1) .. controls (2,0.1) and (1.1,0.4) .. (1.5,0.6);

\draw[thick, latex-] (0.8+0.48*2-0.02,0.01+0.32*1.5+0.04) -- (0.8+0.48*2-0.3,0.01+0.32*1.5+0.6);
\draw[thick,latex-] (0.8+0.48*3.5-0.02,0.01+0.32*2+0.05) -- (0.8+0.48*3.5-0.24,0.01+0.32*2+0.6);
\draw[thick,latex-] (0.8+0.48*5-0.01,0.01+0.32*3+0.05) -- (0.8+0.48*5-0.12,0.01+0.32*3+0.6);
\draw[thick,latex-] (0.8+0.48*6.25-0.01,0.01+0.32*1.5+0.05) -- (0.8+0.48*6.25-0.12,0.01+0.32*1.5+0.6);
\draw[thick,latex-] (0.8+0.48*7-0.005,0.01+0.32*3.5+0.05) -- (0.8+0.48*7-0.06,0.01+0.32*3.5+0.6);
\draw[thick,latex-] (0.8+0.48*8+0.005,0.01+0.32*2+0.05) -- (0.8+0.48*8+0.06,0.01+0.32*2+0.6);
\draw[thick,latex-] (0.8+0.48*9+0.005,0.01+0.32*3.5+0.05) -- (0.8+0.48*9+0.06,0.01+0.32*3.5+0.6);
\draw[thick,latex-] (0.8+0.48*10+0.01,0.01+0.32*2+0.05) -- (0.8+0.48*10+0.12,0.01+0.32*2+0.6);

\node at (6.2,1.1) {$\mathcal{E}_{\eta}$};
\end{scope}

\end{tikzpicture}
\end{center}
\caption{\label{fig5} The bijection $z \mapsto (z, \theta(z,\eta))$ takes $\mathcal{Z}$ to $\mathcal{E}_{\eta}$ and hence $\mathcal{Z}^o$ to $\mathcal{E}_\eta^o$. $\mathcal{E}_\eta^o$ consists of equilibria of ($\mathrm{B_0}$), and these equilibria have real, negative eigenvalues in directions traverse to $\mathcal{E}_\eta$.}
\end{figure}

(ii) ($\mathrm{A_0}$) arises from the differential-algebraic system $\dot z = \Gamma v(z+\alpha\hatt{\beta}^{-1}\hatt{y}),\,\,0=\hatt{\beta}f_{*}(z,\hatt{y},\eta)$ bearing in mind that $\hatt{\beta}$ is nonsingular and that solutions to $\hatt{y}=\theta(z,\eta)$ satisfying $(z,\eta) \in \mathcal{Z} \times (0, \eta_4)$ form a portion of the zero set of $f_{*}(z,\hatt{y},\eta)$ (parts~\ref{pointposzero}~and~\ref{pointfirstbound} above). For each fixed $\eta \in (0, \eta_4)$, ($\mathrm{A_0}$) defines a local flow on $\mathcal{Z}^o$ which includes the periodic orbit $\mathcal{O}_\eta'$ which is nondegenerate (resp., linearly stable) relative to $S_{\mathcal{O}}$ (recall Part~\ref{pointthirdbound} above). The projection of the vector field of ($\mathrm{A_0}$), namely $\Gamma v(z+\alpha\hatt{\beta}^{-1}\theta(z, \eta))$, onto the tangent space of $\mathcal{E}_\eta^o$ is the {\em reduced vector field} associated with this system in the terminology of \cite{Fenichel79}. The reduced vector field has the periodic orbit 
\[
\mathcal{O}_\eta := \{(z,\hatt{y})\colon z \in \mathcal{O}'_\eta, \hatt{y} = \theta(z, \eta)\}\,,
\]
on $\mathcal{E}_\eta^o$. $\mathcal{O}_\eta$ is nondegenerate (resp., linearly stable) relative to $\mathcal{E}_\eta^o$ because $\mathcal{O}'_\eta$ is nondegenerate (resp., linearly stable) relative to $S_{\mathcal{O}}$ and $z \mapsto (z, \theta(z, \eta))$ is a $C^1$ diffeomorphism taking $\mathcal{Z}^o \subseteq S_{\mathcal{O}}$ to $\mathcal{E}_\eta^o$ (see Remark~\ref{remfloqrel}). The situation is illustrated in Figure~\ref{fig6}.

\begin{figure}[h]
\begin{center}
\begin{tikzpicture}
\draw[thick,-latex] (0,0) -- (7,0);
\draw[thick,-latex] (0,0) -- (0,4);
\node at (0,4.2) {$\hatt{y}$};
\draw[thick,-latex] (0,0) -- (3,2);

\draw[thick] (1.5,0.6) .. controls (2.1,0.9) and (3,1.2) .. (4.5,1.2) -- (4.5,1.2) .. controls (5.5,1.2) and (6.6,0.9) .. (6,0.6) -- (6,0.6) .. controls (5.4,0.3) and (5,0.1) ..  (3.5,0.1) -- (3.5,0.1) .. controls (2,0.1) and (1.1,0.4) .. (1.5,0.6);
\fill[black!30] (1.5,0.6) .. controls (2.1,0.9) and (3,1.2) .. (4.5,1.2) -- (4.5,1.2) .. controls (5.5,1.2) and (6.6,0.9) .. (6,0.6) -- (6,0.6) .. controls (5.4,0.3) and (5,0.1) ..  (3.5,0.1) -- (3.5,0.1) .. controls (2,0.1) and (1.1,0.4) .. (1.5,0.6);
\node at (6.2,0.5) {$\mathcal{Z}$};

\draw[-latex, dashed] (6.2,0.8) -- (6.2,2.8);
\draw[-latex, dashed] (1.4,0.5) -- (1.4,2.5);
\node at (7.4,1.8) {$z\mapsto (z,\theta(z,\eta))$};
\draw[thick, ->] (2.3,0.6) .. controls (2.7,0.8) and (3.5,1) .. (4.3,1);
\draw[thick, ->] (5.3,0.6) .. controls (4.9,0.4) and (4.3,0.3) .. (3.3,0.3);
\draw[thick] (2.3,0.6) .. controls (2.7,0.8) and (3.5,1) .. (4.3,1) -- (4.3,1) .. controls (5,1) and (5.7,0.8) .. (5.3,0.6) -- (5.3,0.6) .. controls (4.9,0.4) and (4.3,0.3) .. (3.3,0.3) -- (3.3,0.3) .. controls (2.3,0.3) and (1.9,0.4) .. (2.3,0.6);

\node at (4.7,0.7) {$\mathcal{O}_{\eta}'$};

\begin{scope}[yshift=2cm]
\shade (1.5,0.6) .. controls (2.1,0.9) and (3,1.4) .. (4.5,1.4) -- (4.5,1.4) .. controls (5.5,1.4) and (6.6,0.9) .. (6,0.6) -- (6,0.6) .. controls (5.4,0.3) and (5,0.1) ..  (3.5,0.1) -- (3.5,0.1) .. controls (2,0.1) and (1.1,0.4) .. (1.5,0.6);
\draw[black!20] (1.5,0.6) .. controls (2.1,0.9) and (3,1.4) .. (4.5,1.4) -- (4.5,1.4) .. controls (5.5,1.4) and (6.6,0.9) .. (6,0.6) -- (6,0.6) .. controls (5.4,0.3) and (5,0.1) ..  (3.5,0.1) -- (3.5,0.1) .. controls (2,0.1) and (1.1,0.4) .. (1.5,0.6);
\end{scope}
\begin{scope}[yshift=2.2cm]
\draw[thick, ->] (2.3,0.6) .. controls (2.7,0.8) and (3.5,1) .. (4.3,1);
\draw[thick, ->] (5.3,0.6) .. controls (4.9,0.4) and (4.3,0.3) .. (3.3,0.3);
\draw[thick] (2.3,0.6) .. controls (2.7,0.8) and (3.5,1) .. (4.3,1) -- (4.3,1) .. controls (5,1) and (5.7,0.8) .. (5.3,0.6) -- (5.3,0.6) .. controls (4.9,0.4) and (4.3,0.3) .. (3.3,0.3) -- (3.3,0.3) .. controls (2.3,0.3) and (1.9,0.4) .. (2.3,0.6);

\node at (6.2,0.9) {$\mathcal{E}_{\eta}$};
\node at (4.7,0.7) {$\mathcal{O}_{\eta}$};

\end{scope}

\end{tikzpicture}
\end{center}
\caption{\label{fig6} The bijection $z \mapsto (z, \theta(z,\eta))$, which takes $\mathcal{Z}$ to $\mathcal{E}_{\eta}$, lifts the periodic orbit $\mathcal{O}_{\eta}'$ to $\mathcal{O}_{\eta}$.}
\end{figure}
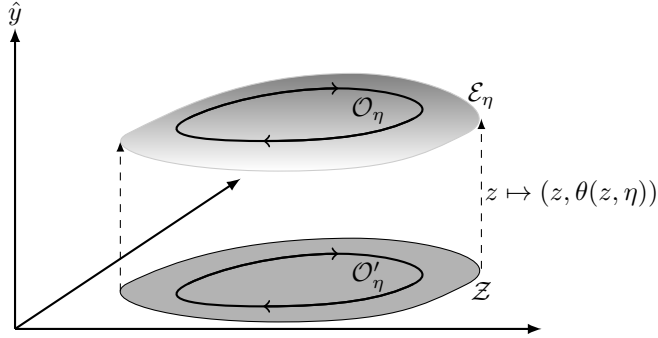

According to Theorem 13.1 in \cite{Fenichel79}, (i) and (ii) together tell us that the periodic orbit $\mathcal{O}_\eta$ ``survives'' perturbation namely, given any $\zeta>0$, we can choose $\epsilon_0>0$ such that for $\epsilon \in [0, \epsilon_0)$, ($\mathrm{A}_\epsilon$) has a periodic orbit $\mathcal{O}_{\eta, \epsilon}$ satisfying $d_H(\mathcal{O}_\eta, \mathcal{O}_{\eta, \epsilon})<\zeta$. In particular, we choose $\epsilon_0>0$ such that $\epsilon \in [0,\epsilon_0)$ implies that $d_H(\mathcal{O}_\eta, \mathcal{O}_{\eta, \epsilon})< d_H(\mathcal{O}_\eta, \partial\mathcal{Z}_+)$, and hence $\mathcal{O}_{\eta, \epsilon} \subseteq \mathcal{Z}_+^o$. 

Moreover, according to Theorem 13.2 in \cite{Fenichel79}, as $\mathcal{O}_\eta$ is nondegenerate (resp., linearly stable) relative to $S$, $\epsilon_1\in (0,\epsilon_0]$ can be chosen such that for all $\epsilon \in [0, \epsilon_1)$, $\mathcal{O}_{\eta, \epsilon}$ is nondegenerate (resp., linearly stable) relative to $S$ as a periodic orbit of ($\mathrm{A}_\epsilon$). The claim about linear stability follows because the nontrivial eigenvalues of ($\mathrm{B_0}$) relative to $S$ at points of $\mathcal{E}_\eta$ all have negative real parts (Part~\ref{pointfourthbound} above).

We now fix $\epsilon \in (0,\epsilon_1)$ and return to $(x,y)$-space. Let $\mathcal{O}_{\eta, \epsilon}^{*} := \phi_0^{-1}\circ \phi_1^{-1}(\mathcal{O}_{\eta, \epsilon})$. Explicitly, 
\[
\mathcal{O}_{\eta, \epsilon}^{*} = \{(x,(\hatt{y},\doublehat{y}))\colon (z,\hatt{y}) \in \mathcal{O}_{\eta, \epsilon},\,\,\doublehat{y} = \mathbf{1}-\delta^{\mathrm{t}}\hatt{y},\,\, x=z+\alpha\hatt{\beta}^{-1}\hatt{y}\}\,.
\]
Since $\mathcal{O}_{\eta, \epsilon} \subseteq \mathcal{Z}_+$, it follows that $\mathcal{O}_{\eta, \epsilon}^{*} \subseteq \mathcal{Z}_+^* \subseteq S_+^*$ (recall (\ref{eqmapstopos1})), namely, $\mathcal{O}_{\eta, \epsilon}^{*}$ is a positive periodic orbit on $S^*$ for the enlarged CRN $\mathcal{R}'$ governed by (\ref{perteq}).

On the other hand, nondegeneracy (resp., linear stability) of $\mathcal{O}_{\eta, \epsilon}$ for ($\mathrm{A}_\epsilon$) relative to $S$, is equivalent (since $\phi_0^{-1}\circ \phi_1^{-1}$ is a $C^1$ diffeomorphism taking $S$ to $S^*$ -- recall (\ref{eqsodash}) and Lemma~\ref{lemfloqdiff}) to nondegeneracy (resp., linear stability) of $\mathcal{O}_{\eta, \epsilon}^{*}$ relative to $S^*$. 

We have thus constructed a family of nondegenerate (resp., linearly stable) periodic orbits of $\mathcal{R}'$ as desired. This completes the proof.

\end{enumerate}

\subsection{Remarks and subproofs}

\begin{remark}
\label{remlocalcoords}
In previous work \cite{banajiCRNosci}, we introduced local coordinates on $S_{\mathcal{O}}$ as follows. Let $r$ be the rank of $\Gamma$, let $\Gamma_0$ be a matrix whose columns form a basis of $\mathrm{im}\,\Gamma$, and define the $r \times r_0$ matrix $Q$ with rank $r$ by $\Gamma = \Gamma_0Q$. Choose an arbitrary point $x_0 \in S_{\mathcal{O}}$ and define $h\colon \mathbb{R}^r \to S_{\mathcal{O}}$ by $h(w)=x_0+\Gamma_0w$. Then $w \in \mathbb{R}^r$ defines a local coordinate on $S_{\mathcal{O}}$ which evolves according to the differential equation
\begin{equation}
\label{eqbasic1}
\dot w = Qv(x_0+\Gamma_0 w)\,.
\end{equation}
We could work with (\ref{eqbasic1}) rather than (\ref{eqbasic0}) as our starting point. However, we avoid explicitly introducing local coordinates as it is not strictly necessary here and \textcolor{black}{introduces an extra transformation for us to track, obscuring} the fundamental geometrical meaning of calculations. Instead, the reader may find it helpful to bear in mind that we can identify $S_{\mathcal{O}}$ with $\mathbb{R}^r$, $\mathcal{Z}$ with a compact subset of $\mathbb{R}^r$, and $\mathcal{Z}^o$ with an open subset of $\mathbb{R}^r$, all via some linear bijection such as $h$.
\end{remark}

\begin{subproof}
\label{phi0cosets}
Given any $z_0 \in \mathbb{R}^n$, $y_0 \in \mathbb{R}^{m+k}$, $t\in \mathbb{R}^{r_0}$ and $s \in \mathbb{R}^m$, 
\begin{eqnarray*}
\phi_0^{-1}\left[\left(\begin{array}{c}z_0\\y_0\end{array}\right) + \left(\begin{array}{cc}\Gamma&0\\0&\beta\end{array}\right)\left(\begin{array}{c}t\\s\end{array}\right)\right]
&=& \left(\begin{array}{c}z_0+\Gamma t + \alpha\hatt{\beta}^{-1}(\hatt{y}_0+\hatt{\beta} s)\\y_0+\beta s\end{array}\right)\\
&=& \left(\begin{array}{c}z_0+\alpha\hatt{\beta}^{-1}\hatt{y}_0\\y_0\end{array}\right) + \left(\begin{array}{cc}\Gamma&\alpha\\0&\beta\end{array}\right)\left(\begin{array}{c}t\\s\end{array}\right)\\
&=& \phi_0^{-1}\left(\begin{array}{c}z_0\\y_0\end{array}\right) + \left(\begin{array}{cc}\Gamma&\alpha\\0&\beta\end{array}\right)\left(\begin{array}{c}t\\s\end{array}\right)
\end{eqnarray*}
\hfill$\square$
\end{subproof}

\begin{subproof}
\label{yhathat}
The claim that $\delta^{\mathrm{t}}\hatt{y}+\doublehat{y}$ is constant along trajectories of (\ref{eqzyfull}) is vacuously true in the case $k=0$. So, suppose that $k>0$. Let $P:=[-\doublehat{\beta}\hatt{\beta}^{-1}|I_k] \in \mathbb{R}^{k \times (m+k)}$. By a quick calculation, $P\beta=0$, and since $P$ has rank $k$ and $\mathrm{dim}(\mathrm{ker}\,\beta^{\mathrm{t}})=k$ by the rank-nullity theorem, the rows of $P$ must form a basis of $\mathrm{ker}\,\beta^{\mathrm{t}}$. Multiplying $\dot y = \beta q(x,y)$ by $P$ on the left gives $P\dot y = 0$, namely $Py = -\doublehat{\beta}\hatt{\beta}^{-1}\hatt{y} + \doublehat{y}$ is constant along trajectories. Thus, the value of $\doublehat{y}$ along a trajectory at any point is specified by the value of $\hatt{y}$ on the trajectory at that point, and an additional parameter, $\delta^{\mathrm{t}}\hatt{y} + \doublehat{y} \in \mathbb{R}^k$. \hfill$\square$
\end{subproof}

\begin{subproof}
\label{phi0phi1}
Let $x_0$ be some arbitrary point on $S_{\mathcal{O}}$. Note that
\[
\phi_0^{-1}\circ \phi_1^{-1}(S) = \{(x,y)\colon x-\alpha\hatt{\beta}^{-1}\hatt{y} \in S_{\mathcal{O}},\, \delta^{\mathrm{t}}\hatt{y}+\doublehat{y}=\mathbf{1}\}\,.
\]
Thus $(x,y) \in \phi_0^{-1}\circ \phi_1^{-1}(S)$, implies that $x-\alpha\hatt{\beta}^{-1}\hatt{y} = x_0 + \Gamma t$ for some $t \in \mathbb{R}^{r_0}$, namely $x=x_0+\Gamma t + \alpha\hatt{\beta}^{-1}\hatt{y}$, and $\doublehat{y} = \mathbf{1}-\delta^{\mathrm{t}}\hatt{y}$, so that
\[
\left(\begin{array}{c}x\\\hatt{y}\\\doublehat{y}\end{array}\right) = \left(\begin{array}{c}x_0\\0\\\mathbf{1}\end{array}\right) + \left(\begin{array}{cc}\Gamma&\alpha\\0&\hatt{\beta}\\0&\doublehat{\beta}\end{array}\right)\left(\begin{array}{c}t\\\hatt{\beta}^{-1}\hatt{y}\end{array}\right) \in \left(\begin{array}{c}x\\0\\\mathbf{1}\end{array}\right) + \mathrm{im}\,\Gamma'\,.
\]
Thus $(x,y) \in \phi_0^{-1}\circ \phi_1^{-1}(S) \Rightarrow (x,y) \in S^*$. On the other hand, $(x,y) \in S^*$ means, by definition, that
\[
\left(\begin{array}{c}x\\\hatt{y}\\\doublehat{y}\end{array}\right) = \left(\begin{array}{c}x_0\\0\\\mathbf{1}\end{array}\right) + \left(\begin{array}{cc}\Gamma&\alpha\\0&\hatt{\beta}\\0&\doublehat{\beta}\end{array}\right)\left(\begin{array}{c}t\\s\end{array}\right)\,,
\]
for some $t \in \mathbb{R}^{r_0}$ and $s \in \mathbb{R}^m$, from which we see that $x-\alpha\hatt{\beta}^{-1}\hatt{y} = x_0 + \Gamma t$ (namely, $x-\alpha\hatt{\beta}^{-1}\hatt{y} \in S_{\mathcal{O}}$), and $\delta^{\mathrm{t}}\hatt{y}+\doublehat{y}=\mathbf{1}$. Thus $(x,y) \in S^* \Rightarrow (x,y) \in \phi_0^{-1}\circ \phi_1^{-1}(S)$, confirming the claim. \hfill$\square$
\end{subproof}

\begin{subproof}
\label{zplus}
Let $D = d_H(\mathcal{Z}, \partial \mathbb{R}^n_{\geq 0})$, and let $\|\cdot\|$ refer to the matrix norm induced by the Euclidean norm. If $k=0$, then define $y_{pos} := D/(2\|\alpha\hatt{\beta}^{-1}\|)\,.$ Otherwise, set 
\[
y_{pos} := \min\left\{\frac{D}{2\|\alpha\hatt{\beta}^{-1}\|}, \,\, \frac{1}{2\|\delta^{\mathrm{t}}\|}\right\}\,.
\]
Then $(z,\hatt{y}) \in \mathcal{Z}_+$ implies that (i) $z+\alpha\hatt{\beta}^{-1}\hatt{y} \gg 0$ and (ii) $\delta^{\mathrm{t}}\hatt{y} \ll \mathbf{1}$ and, consequently, $\mathbf{1}-\delta^{\mathrm{t}}\hatt{y} \gg 0\,.$ Consequently,
\[
\phi_0^{-1}\circ \phi_1^{-1}(\mathcal{Z}_+)=\{(z+\alpha\hatt{\beta}^{-1}\hatt{y},(\hatt{y},\mathbf{1}-\delta^{\mathrm{t}}\hatt{y}))\colon z \in \mathcal{Z},\,\,\hatt{y} \gg 0,\,\, |\hatt{y}| \leq y_{pos}\} \subseteq \mathbb{R}^n_{\gg 0} \times \mathbb{R}^{m+k}_{\gg 0}\,.
\]
\hfill$\square$
\end{subproof}

\begin{subproof}
\label{ftog}
Solving $f(x,y, \eta) = 0$ with the assumption that $x \gg 0$ and $y \gg 0$ gives $y^{\beta^{\mathrm{t}}} = \bm{\eta}^{\hatt{\beta}^{\mathrm{t}}} \circ x^{-\alpha^{\mathrm{t}}}$. This can be written $\hatt{y}^{\hatt{\beta}^{\mathrm{t}}} = \bm{\eta}^{\hatt{\beta}^{\mathrm{t}}} \circ x^{-\alpha^{\mathrm{t}}}\circ \doublehat{y}^{-\doublehat{\beta}^{\mathrm{t}}}$. Taking logs gives:
\[
\hatt{\beta}^{\mathrm{t}}\ln \hatt{y} = \hatt{\beta}^{\mathrm{t}}\ln(\bm{\eta}) -\alpha^{\mathrm{t}} \ln x -\doublehat{\beta}^{\mathrm{t}} \ln \doublehat{y}\,.
\]
Multiplying through by $(\hatt{\beta}^{\mathrm{t}})^{-1}$ and exponentiating again gives the result. \hfill$\square$
\end{subproof}

\begin{subproof}
\label{theta}
First, we calculate that $g_{*}(z,0,0)=0$ for any $z$ in the domain of $g_{*}$, and that
\[
D_{\hatt{y}}g_{*}(z,0,0) = \mathrm{id}\,.
\]
We now apply the extension of the IFT in Lemma~\ref{lemIFT1}. Note that given any neighbourhood $V$ of $\mathcal{Z} \times \{0\}$ in $\mathbb{R}^n \times \mathbb{R}$, there exists $\eta_1>0$ such that $\mathcal{Z} \times (-\eta_1,\eta_1) \subseteq V$. Consequently, since $g_{*}$ is $C^2$, according to Lemma~\ref{lemIFT1} there exists $\eta_1>0$, $y_{max} \in (0, y_{pos}]$ and a $C^2$ function $\theta\colon \mathcal{Z} \times (-\eta_1,\eta_1) \to B_{y_{max}}$ whose graph is precisely the zero-set of $g_{*}$ in $\mathcal{Z} \times B_{y_{max}} \times (-\eta_1,\eta_1)$, namely,
\[
\{(z,\hatt{y},\eta) \in \mathcal{Z} \times B_{y_{max}} \times (-\eta_1,\eta_1)\,\colon\, g_{*}(z,\hatt{y},\eta)=0\} = \{(z,\theta(z,\eta),\eta)\,\colon\, (z,\eta) \in \mathcal{Z} \times (-\eta_1,\eta_1)\}\,.
\] \hfill$\square$
\end{subproof}

\begin{subproof}
\label{ypos}
We find $\eta_2\in (0,\eta_1]$, such that $\theta(z,\eta) \gg 0$ for $(z,\eta) \in \mathcal{Z} \times (0,\eta_2)$. Consider the function $\hatt{\theta}\colon \mathcal{Z} \times (-\eta_1,\eta_1)$ defined by
\[
\hatt{\theta}(z,\eta)=\left\{\begin{array}{rcl}\theta(z,\eta)/\eta & (\eta\neq 0),\\D_{\eta}\theta(z,0) & (\eta=0).\end{array}\right.
\]
Since $\theta$ is $C^2$ on its domain and $\theta(z, 0)=0$ for all $z \in \mathcal{Z}$ we have, by Lemma~\ref{lemuniform}, that $\hatt{\theta}$ is continuous on its domain. On the other hand $\theta(\mathcal{Z},0) \subseteq \mathbb{R}^m_{\gg 0}$: implicitly differentiating $g_{*}(z,\theta(z, \eta),\eta)=0$ w.r.t. $\eta$, evaluating at $(z,0)$ for any $z \in \mathcal{Z}$, recalling that $D_{\hatt{y}}g_{*}(z,0,0) = \mathrm{id}$, and using the fact that $\theta(z,0)=0$ gives
\[
\hatt{\theta}(z,0) = D_{\eta}\theta(z,0) = -(D_yg_{*}(z,0,0))^{-1}D_{\eta}g_{*}(z,0,0) = -D_{\eta}g_{*}(z,0,0) = z^{\gamma} \gg 0\,.
\]
If $\hatt{\theta}(z,\eta)\gg 0$ on $\mathcal{Z} \times (-\eta_1,\eta_1)$, then $\theta(z,\eta)\gg 0$ on $\mathcal{Z} \times (0,\eta_1)$ and we define $\eta_2=\eta_1$. Otherwise, we define
\[
\eta_2:=\inf\{\eta \in (0,\eta_1) \colon (z,\eta) \in \theta^{-1}(\partial\mathbb{R}^m_{\geq 0})\} = \inf\{\eta \in (0,\eta_1) \colon (z,\eta) \in \hatt{\theta}^{-1}(\partial\mathbb{R}^m_{\geq 0})\} > 0\,.
\]
The last inequality follows by continuity of $\hatt{\theta}$: the pre-image of the closed set $\partial \mathbb{R}^m_{\geq 0}$ under $\hatt{\theta}$ is closed, and would otherwise have to include a point in $\mathcal{Z} \times \{0\}$, contradicting $\hatt{\theta}(\mathcal{Z},0)\subseteq \mathbb{R}^m_{\gg 0}$. \hfill$\square$
\end{subproof}

\begin{subproof}
\label{normhyp}
We can rewrite $f$ as a function of four variables $x$, $\hatt{y}, \doublehat{y}$ and $\eta$, namely
\[
\overline{f}(x,\hatt{y}, \doublehat{y},\eta) := \bm{\eta}^{-\hatt{b}^{\mathrm{t}}}\circ x^{a^{\mathrm{t}}}\circ \hatt{y}^{\hatt{b}^{\mathrm{t}}}\circ \doublehat{y}^{\doublehat{b}^{\mathrm{t}}} - \bm{\eta}^{-{\hatt{b}'}^{\mathrm{t}}}\circ x^{{a'}^{\mathrm{t}}}\circ \hatt{y}^{{\hatt{b}'}^{\mathrm{t}}}\circ \doublehat{y}^{{\doublehat{b}'}^{\mathrm{t}}}\,.
\]
Noting that $f_{*}(z, \hatt{y}, \eta) = \overline{f}(z+\alpha\hatt{\beta}^{-1}\hatt{y},\hatt{y},\mathbf{1}-\delta^{\mathrm{t}}\hatt{y}, \eta)$, we get
\begin{eqnarray*}
W(z,\eta)&=&\hatt{\beta} D_1 \overline{f}(z+\alpha\hatt{\beta}^{-1}\theta(z,\eta),\theta(z,\eta),\mathbf{1}-\delta^{\mathrm{t}}\theta(z,\eta), \eta)\alpha\hatt{\beta}^{-1}\\
&&\hspace{1cm}+\hatt{\beta} D_2 \overline{f}(z+\alpha\hatt{\beta}^{-1}\theta(z,\eta),\theta(z,\eta),\mathbf{1}-\delta^{\mathrm{t}}\theta(z,\eta),\eta)\\
&&\hspace{2cm}-\hatt{\beta} D_3 \overline{f}(z+\alpha\hatt{\beta}^{-1}\theta(z,\eta),\theta(z,\eta),\mathbf{1}-\delta^{\mathrm{t}}\theta(z,\eta),\eta)\delta^{\mathrm{t}}\,.
\end{eqnarray*}
$W\colon \mathcal{Z} \times (0, \eta_3) \to \mathbb{R}^{m \times m}$ is defined and continuous. We want to show that $\eta W$ has a continuous extension to $\mathcal{Z} \times [0,\eta_3)$, namely, there exists a continuous function $\overline{W}\colon \mathcal{Z} \times [0,\eta_3) \to \mathbb{R}^{m \times m}$ which coincides with $\eta W$ on $\mathcal{Z} \times (0,\eta_3)$. Further, $\overline{W}(z,0)$ has real, negative eigenvalues for each fixed $z \in \mathcal{Z}$. By continuity of $\overline{W}$ and compactness of $\mathcal{Z}$, it will follow that there exists $\eta_4 \in (0, \eta_3]$ such that for $(z,\eta)\in \mathcal{Z} \times (0, \eta_4)$, eigenvalues of $\overline{W}(z, \eta)$ lie in the open left half of the complex plane. The same must hold for $W(z,\eta)$ whose eigenvalues are positive multiples of those of $\overline{W}(z,\eta)$. \\\\
Note first that $f_{*}(z,\theta(z,\eta),\eta)=0$ for $(z,\eta) \in \mathcal{Z} \times (0,\eta_3)$, implying that
\begin{eqnarray*}
T'(z,\eta)& :=& \bm{\eta}^{-\hatt{b}^{\mathrm{t}}}\circ (z+\alpha\hatt{\beta}^{-1}\theta(z,\eta))^{a^{\mathrm{t}}}\circ \theta(z, \eta)^{\hatt{b}^{\mathrm{t}}}\circ (\mathbf{1}-\delta^{\mathrm{t}}\theta(z, \eta))^{\doublehat{b}^{\mathrm{t}}} \\&=& \bm{\eta}^{-{\hatt{b}'}^{\mathrm{t}}}\circ (z+\alpha\hatt{\beta}^{-1}\theta(z,\eta))^{{a'}^{\mathrm{t}}}\circ \theta(z, \eta)^{{\hatt{b}'}^{\mathrm{t}}}\circ (\mathbf{1}-\delta^{\mathrm{t}}\theta(z, \eta))^{{\doublehat{b}'}^{\mathrm{t}}}\,.
\end{eqnarray*}
$T'\colon \mathcal{Z} \times (0,\eta_3) \to \mathbb{R}^m$ has a continuous extension to $\mathcal{Z} \times [0,\eta_3)$. To see this note, via Lemma~\ref{lemuniform} and Subproof~\ref{ypos}, that 
\[
\theta(z, \eta) = \eta(z^\gamma + R(z, \eta))
\]
where $R(z, \eta)$ is continuous on $\mathcal{Z} \times [0, \eta_3)$ and satisfies $R(z, 0) = 0$. Moreover $\hatt{\theta}(z,\eta):=z^\gamma + R(z, \eta)$ is continuous and positive on $\mathcal{Z} \times [0, \eta_3)$ (see Subproof~\ref{ypos}). So 
\[
T(z,\eta) := (z+\alpha\hatt{\beta}^{-1}\theta(z, \eta))^{a^{\mathrm{t}}}\circ \hatt{\theta}(z,\eta)^{\hatt{b}^{\mathrm{t}}}\circ (\mathbf{1}-\delta^{\mathrm{t}}\theta(z, \eta))^{\doublehat{b}^{\mathrm{t}}}
\]
is defined and continuous on $\mathcal{Z} \times [0, \eta_3)$ and clearly extends $T'(z,\eta)$. Simple evaluation gives 
\[
T(z,0)= z^{a^{\mathrm{t}}}\circ z^{\hatt{b}^{\mathrm{t}}\gamma} \gg 0\,.
\]
We are now ready to compute $\overline{W}(z,\eta)$. Differentiating $\overline{f}$ w.r.t. its first argument (see Example~\ref{exdiffmon}) gives
\[
D_1\overline{f}(x,\hatt{y},\doublehat{y},\eta) = \mathrm{diag}(T_{1})a^{\mathrm{t}}\mathrm{diag}(\mathbf{1}/x) - \mathrm{diag(T_2})a'^{\mathrm{t}}\mathrm{diag}(\mathbf{1}/x)
\]
where $T_1$ and $T_2$ are abbreviations for the first and second terms in $\overline{f}$ respectively. But when $\hatt{y}=\theta(z,\eta)$ and $\doublehat{y}=\mathbf{1}-\delta^{\mathrm{t}}\theta(z,\eta)$, then $\overline{f} = 0$, and $T_1=T_2=T(z,\eta)$. Define $D_{z,\eta}:= \mathrm{diag}(T(z,\eta))$ on $\mathcal{Z} \times [0,\eta_3)$. We calculate:
\[
\eta D_1 \overline{f}(z+\alpha\hatt{\beta}^{-1}\theta(z,\eta),\theta(z,\eta),\mathbf{1}-\delta^{\mathrm{t}}\theta(z,\eta), \eta) = -\eta D_{z,\eta}\alpha^{\mathrm{t}}\mathrm{diag}(\mathbf{1}/(z+\alpha\hatt{\beta}^{-1}\theta(z, \eta)))\,.
\]
This quantity is defined and continuous on $\mathcal{Z} \times [0,\eta_3)$ since $(z,\eta) \in \mathcal{Z} \times [0,\eta_3)$ implies that $(z,\theta(z,\eta)) \in \mathcal{Z}_+$ and consequently $z+\alpha\hatt{\beta}^{-1}\theta(z, \eta) \gg 0$ (see Subproof~\ref{zplus}). Similarly,
\[
\eta D_2 \overline{f}(z+\alpha\hatt{\beta}^{-1}\theta(z,\eta),\theta(z,\eta),\mathbf{1}-\delta^{\mathrm{t}}\theta(z,\eta),\eta) = -D_{z,\eta}\hatt{\beta}^{\mathrm{t}}\mathrm{diag}(\mathbf{1}/\hatt{\theta}(z,\eta))\,.
\]
This is again defined and continuous on $\mathcal{Z} \times [0,\eta_3)$ (we have already observed in Subproof~\ref{ypos} that $\hatt{\theta}$ is continuous and positive on $\mathcal{Z} \times [0,\eta_3)$.) Finally,
\[
\eta D_3 \overline{f}(z+\alpha\hatt{\beta}^{-1}\theta(z,\eta),\theta(z,\eta),\mathbf{1}-\delta^{\mathrm{t}}\theta(z,\eta),\eta) = -\eta D_{z,\eta}\doublehat{\beta}^{\mathrm{t}}\mathrm{diag}(\mathbf{1}/(\mathbf{1}-\delta^{\mathrm{t}}\theta(z,\eta)))\,.
\]
(Note that if $k=0$, then according to our conventions $D_3\overline{f}$ is a zero matrix.) Again, this quantity is defined and continuous on $\mathcal{Z} \times [0,\eta_3)$ since $(z,\eta) \in \mathcal{Z} \times [0,\eta_3)$ implies that $(z,\theta(z,\eta)) \in \mathcal{Z}_+$ and consequently $\mathbf{1}-\delta^{\mathrm{t}}\theta(z,\eta)$ is positive (see Subproof~\ref{zplus}). We thus have that
\begin{eqnarray*}
\overline{W}(z,\eta) &:=& -\eta\hatt{\beta}D_{z,\eta}\alpha^{\mathrm{t}}\mathrm{diag}(\mathbf{1}/(z+\alpha\hatt{\beta}^{-1}\theta(z, \eta)))\alpha\hatt{\beta}^{-1} \,\,-\,\, \hatt{\beta} D_{z,\eta}\hatt{\beta}^{\mathrm{t}}\mathrm{diag}(\mathbf{1}/\hatt{\theta}(z,\eta))\\
&&\hspace{4cm}+\,\,\eta\hatt{\beta}D_{z,\eta}\doublehat{\beta}^{\mathrm{t}}\mathrm{diag}(\mathbf{1}/(\mathbf{1}-\delta^{\mathrm{t}}\theta(z,\eta)))\delta^{\mathrm{t}}\,
\end{eqnarray*}
is defined and continous on $\mathcal{Z} \times [0,\eta_3)$ and coincides with $\eta W(z,\eta)$ on $(0, \eta_3)$. Moreover, as calculated in Subproof~\ref{ypos}, $\hatt{\theta}(z,0)=z^\gamma\gg 0$, and so
\[
\overline{W}(z,0)= -\hatt{\beta} D_{z,0}\hatt{\beta}^{\mathrm{t}}\mathrm{diag}(\mathbf{1}/z^\gamma)\,.
\]
Defining $Y_z:=\mathrm{diag}(\mathbf{1}/z^\gamma)$, we have
\[
\overline{W}(z,0)= -Y_z^{-1/2}[Y_z^{1/2}\hatt{\beta} D_{z,0}\hatt{\beta}^{\mathrm{t}}Y_z^{1/2}]Y_z^{1/2} = -Y_z^{-1/2}(M_zM_z^{\mathrm{t}})Y_z^{1/2}\,,
\]
where $M_z:=Y_z^{1/2}\hatt{\beta}(D_{z,0})^{1/2}$ is nonsingular as each of its factors is nonsingular. From this we see that $\overline{W}(z,0)$ is similar to a negative definite matrix, and so its eigenvalues are real and negative as claimed. \hfill$\square$
\end{subproof}

\section{Elucidating the proof via an example}
\label{seceluc}

\textcolor{black}{We now return to the example in Section~\ref{secexample} and use it to elucidate the proof of Theorem~\ref{mainthm}. We derive for this example many of the matrices, sets and functions which occur in the proof above. Recall the CRNs
\begin{equation}
X+Z \overset{\scriptstyle{k_1}}\longrightarrow 2Y \overset{\scriptstyle{k_2}}\longrightarrow X+Y, \quad 0 \overset{\scriptstyle{k_3}}{\underset{\scriptstyle{k_4}}\rightleftharpoons} X, \quad 0 \overset{\scriptstyle{k_5}}{\underset{\scriptstyle{k_6}}\rightleftharpoons} Y, \quad 0 \overset{\scriptstyle{k_7}}{\underset{\scriptstyle{k_8}}\rightleftharpoons} Z\,, \tag{\mbox{$\mathcal{R}_1$}}
\end{equation}
and
\begin{equation}
X+Z \overset{\scriptstyle{k_1}}\longrightarrow 2Y \overset{\scriptstyle{k_2}}\longrightarrow X+Y, \,\,\, 0 \overset{\scriptstyle{k_3}}{\underset{\scriptstyle{k_4}}\rightleftharpoons} X, \,\,\, 0 \overset{\scriptstyle{k_5}}{\underset{\scriptstyle{k_6}}\rightleftharpoons} Y, \,\,\, 0 \overset{\scriptstyle{k_7}}{\underset{\scriptstyle{k_8}}\rightleftharpoons} Z, \,\,\, Y \overset{\scriptstyle{k_9}}{\underset{\scriptstyle{k_{10}}}\rightleftharpoons} U+V, \,\,\, U+X \overset{\scriptstyle{k_{11}}}{\underset{\scriptstyle{k_{12}}}\rightleftharpoons} 2V+W\,. \tag{\mbox{$\mathcal{R}_2$}}
\end{equation}
We fix the order $[X,Y,Z,U,V,W]$ on the chemical species and let the corresponding lower case letters refer to chemical concentrations. The stoichiometric matrices (Part~\ref{pointbasic} above) of $\mathcal{R}_1$ and $\mathcal{R}_2$ are, respectively
\[
\Gamma = \left(\begin{array}{rrrrr}-1&1&1&0&0\\2&-1&0&1&0\\-1&0&0&0&1\end{array}\right) \quad \mbox{and}\quad\Gamma' = \left(\begin{array}{rrrrrrr}-1&1&1&0&0&0&-1\\2&-1&0&1&0&-1&0\\-1&0&0&0&1&0&0\\0&0&0&0&0&1&-1\\0&0&0&0&0&1&2\\0&0&0&0&0&0&1\end{array}\right)\,.
\]
We also have following matrices defined in Part~\ref{pointbasic} of the proof:
\[
a=\left(\begin{array}{cc}0&1\\1&0\\0&0\end{array}\right), \,\, a'=\left(\begin{array}{cc}0&0\\0&0\\0&0\end{array}\right),\,\,\alpha = a'-a=\left(\begin{array}{rr}0&-1\\-1&0\\0&0\end{array}\right)\,,
\]
and
\[
b=\left(\begin{array}{cc}0&1\\0&0\\0&0\end{array}\right), \,\, b'=\left(\begin{array}{cc}1&0\\1&2\\0&1\end{array}\right),\,\,\beta = b'-b=\left(\begin{array}{rr}1&-1\\1&2\\0&1\end{array}\right),\,\,\hatt{\beta}=\left(\begin{array}{rr}1&-1\\1&2\end{array}\right),\,\,\doublehat{\beta}=\left(\begin{array}{rr}0&1\end{array}\right)\,.
\]
From these we obtain (Parts~\ref{pointreduce}~and~\ref{pointposzero} above):
\[
\gamma = -(\alpha\hatt{\beta}^{-1})^\mathrm{t} = \left(\begin{array}{rrr}-1&2&0\\1&1&0\end{array}\right), \quad \mbox{and} \quad \delta=-(\doublehat{\beta}\hatt{\beta}^{-1})^{\mathrm{t}} = \left(\begin{array}{r}1/3\\-1/3\end{array}\right)\,.
\]
The added reactions give rise to the conservation relation (Part~\ref{pointreduce})
\[
u/3-v/3+w=\mathrm{const.}
\]
$v$, the rate vector of $\mathcal{R}_1$, is given by
\[
(k_1xz,\,\, k_2y^2,\,\, k_3-k_4x,\,\, k_5-k_6y,\,\, k_7-k_8z)^\mathrm{t}\,\,.
\]
Thus, the ODE system for $\mathcal{R}_1$ is
\begin{equation}
\label{origCRN}
\left(\begin{array}{c}\dot x\\\dot y\\\dot z\end{array}\right) = \left(\begin{array}{rrrrr}-1&1&1&0&0\\2&-1&0&1&0\\-1&0&0&0&1\end{array}\right)\left(\begin{array}{c}k_1xz\\ k_2y^2\\ k_3-k_4x\\ k_5-k_6y\\ k_7-k_8z\end{array}\right)
\end{equation}
As $\Gamma$ has rank $3$, the stoichiometric subspace is the whole of $(x,y,z)$-space, and consequently, the unique positive stoichiometry class of the system is the positive orthant in $(x,y,z)$-space. We know that for some choice of rate constants this system has a linearly stable periodic orbit. \\\\
Turning to $\mathcal{R}_2$, the rate vector $q$ of the new reactions, set according to Part~\ref{pointrates}, is
\[
\frac{1}{\epsilon}\left(\begin{array}{c}y-\eta^{-2}uv\\\eta^{-1}xu-\eta^{-2}v^2w\end{array}\right)\,.
\]
Thus the ODE system for $\mathcal{R}_2$ is
\begin{equation}
\label{newCRN}
\left(\begin{array}{c}\dot x\\\dot y\\\dot z\\\dot u\\\dot v \\\dot w\end{array}\right) = \left(\begin{array}{rrrrrrr}-1&1&1&0&0&0&-1\\2&-1&0&1&0&-1&0\\-1&0&0&0&1&0&0\\0&0&0&0&0&1&-1\\0&0&0&0&0&1&2\\0&0&0&0&0&0&1\end{array}\right)\left(\begin{array}{c}k_1xz\\ k_2y^2\\ k_3-k_4x\\ k_5-k_6y\\ k_7-k_8z\\(y-\eta^{-2}uv)/\epsilon\\(\eta^{-1}xu-\eta^{-2}v^2w)/\epsilon\end{array}\right)
\end{equation}
Theorem~\ref{mainthm} tells us that for sufficiently small $\epsilon$ and $\eta$ (\ref{newCRN}) has a positive periodic orbit on the set defined by $u/3-v/3+w=1$, and that this periodic orbit is linearly stable relative to its stoichiometry class.\\\\ 
Note at this stage that (\ref{newCRN}) does {\em not} have the form of a singular perturbation problem if we treat $\epsilon$ as a small parameter. In order to remedy this we follow Part~\ref{pointcoords} of the proof and introduce new variables $(x',y',z')$ defined by 
\[
x'=x+v/3-u/3, \,\,y'=y+v/3+2u/3,\,\,z'=z\,.
\]
Conversely
\[
x=x'-v/3+u/3, \,\,y=y'-v/3-2u/3,\,\,z=z'\,.
\]
The dynamical system in the new variables becomes:
\[
\left(\begin{array}{c}\dot{x'}\\\dot{y'}\\\dot{z'}\\\dot{u}\\\dot{v}\\\dot{w}\end{array}\right)\,\,=\,\,
\left(\begin{array}{rrrrrrr}-1&1&1&0&0&0&0\\2&-1&0&1&0&0&0\\-1&0&0&0&1&0&0\\0&0&0&0&0&1&-1\\0&0&0&0&0&1&2\\0&0&0&0&0&0&1\end{array}\right)\left(\begin{array}{c}k_1(x'-v/3+u/3)z'\\ k_2(y'-v/3-2u/3)^2\\ k_3-k_4(x'-v/3+u/3)\\ k_5-k_6(y'-v/3-2u/3)\\ k_7-k_8z'\\((y'-v/3-2u/3)-\eta^{-2}uv)/\epsilon\\(\eta^{-1}(x'-v/3+u/3)u-\eta^{-2}v^2w)/\epsilon\end{array}\right)\,.
\]
Note the block diagonal structure of the new ``stoichiometric matrix''. The system now has the form of a singular perturbation problem with $\epsilon$ as the perturbation parameter. $w$ can be eliminated by fixing the conserved quantity and setting $w = 1-u/3-v/3$ to give:
\begin{equation}
\label{fullCRN}
\begin{array}{rcl}
\left(\begin{array}{c}\dot{x'}\\\dot{y'}\\\dot{z'}\end{array}\right)\,\,&=&\,\,
\left(\begin{array}{rrrrr}-1&1&1&0&0\\2&-1&0&1&0\\-1&0&0&0&1\end{array}\right)\left(\begin{array}{c}k_1(x'-v/3+u/3)z'\\ k_2(y'-v/3-2u/3)^2\\ k_3-k_4(x'-v/3+u/3)\\ k_5-k_6(y'-v/3-2u/3)\\ k_7-k_8z'\end{array}\right)\,,\\
\epsilon\left(\begin{array}{c}\dot{u}\\\dot{v}\end{array}\right)\,\,&=&\,\,
\left(\begin{array}{rr}1&-1\\1&2\\0&1\end{array}\right)\left(\begin{array}{c}(y'-v/3-2u/3)-\eta^{-2}uv\\\eta^{-1}(x'-v/3+u/3)u-\eta^{-2}v^2(1-u/3-v/3)\end{array}\right)\,.
\end{array}
\end{equation}
This singularly perturbed set of ODEs is, for this example, the system referred to as ($\mathrm{A}_\epsilon$) in Part~\ref{pointSPT} of the proof. The function $f_{*}$ referred to in Part~\ref{pointrates} is, in this case, the following function of $x',y',z',u,v$ and $\eta$ which occurs in the second part of (\ref{fullCRN}).
\[
\left(\begin{array}{c}(y'-v/3-2u/3)-\eta^{-2}uv\\\eta^{-1}(x'-v/3+u/3)u-\eta^{-2}v^2(1-u/3-v/3)\end{array}\right)\,.
\]
In the limit $\epsilon \to 0$, (\ref{fullCRN}) reduces to the differential-algebraic system
\begin{equation}
\label{Dalg}
\begin{array}{rcl}
\left(\begin{array}{c}\dot{x'}\\\dot{y'}\\\dot{z'}\end{array}\right)\,\,&=&\,\,
\left(\begin{array}{rrrrr}-1&1&1&0&0\\2&-1&0&1&0\\-1&0&0&0&1\end{array}\right)\left(\begin{array}{c}k_1(x'-v/3+u/3)z'\\ k_2(y'-v/3-2u/3)^2\\ k_3-k_4(x'-v/3+u/3)\\ k_5-k_6(y'-v/3-2u/3)\\ k_7-k_8z'\end{array}\right)\,,\\
0\,\,&=&\,\,
\left(\begin{array}{c}(y'-v/3-2u/3)-\eta^{-2}uv\\\eta^{-1}(x'-v/3+u/3)u-\eta^{-2}v^2(1-u/3-v/3)\end{array}\right)\,,
\end{array}
\end{equation}
while the associated fast-time system gives the ODEs:
\begin{equation}
\label{SingLim}
\begin{array}{rcl}
\left(\begin{array}{c}\dot{x'}\\\dot{y'}\\\dot{z'}\end{array}\right)\,\,&=&\,\,0,\\
\left(\begin{array}{c}\dot{u}\\\dot{v}\end{array}\right)\,\,&=&\,\,
\left(\begin{array}{rr}1&-1\\1&2\\0&1\end{array}\right)\left(\begin{array}{c}(y'-v/3-2u/3)-\eta^{-2}uv\\\eta^{-1}(x'-v/3+u/3)u-\eta^{-2}v^2(1-u/3-v/3)\end{array}\right)\,.
\end{array}
\end{equation}
(\ref{SingLim}) is the system referred to as ($\mathrm{B}_0$) in Part~\ref{pointSPT} of the proof.\\\\
When we observe (\ref{Dalg}), our instinct is to ``solve'' the algebraic equations for $u$ and $v$ in terms of $x', y'$ and $z'$, and substitute into the differential equations. Several steps of the proof are focussed on showing that for fixed small $\eta$ this is indeed possible. Note, however, that even for this simple example explicit solution is difficult, and the sensible approach is to appeal to the implicit function theorem to solve the equations locally.
\\\\
The function $g$ defined in Part~\ref{pointposzero} whose positive zeros are precisely the positive equilibria of the new reactions takes the form
\[
g(x,y,z,u,v,w,\eta)=\left(\begin{array}{c}u\\v\end{array}\right) - \eta\left(\begin{array}{c}x^{-1/3}y^{2/3}w^{1/3}\\x^{1/3}y^{1/3}w^{-1/3}\end{array}\right)\,.
\]
Substituting in $w=1-u/3+v/3$ from the new conservation relation and using the new variables $x', y', z'$ we obtain for this example the function referred to as $g_*$ in Part~\ref{pointposzero} of the proof:
\[
g_*(x',y',z',u,v,\eta) = \left(\begin{array}{c}u\\v\end{array}\right) - \eta\left(\begin{array}{c}(x'-\frac{v}{3}+\frac{u}{3})^{-1/3}(y'-\frac{v}{3}-2\frac{u}{3})^{2/3}(1-\frac{u}{3}+\frac{v}{3})^{1/3}\\(x'-\frac{v}{3}+\frac{u}{3})^{1/3}(y'-\frac{v}{3}-2\frac{u}{3})^{1/3}(1-\frac{u}{3}+\frac{v}{3})^{-1/3}\end{array}\right)\,\,.
\]
Note that $g_*$ is defined provided $x'-v/3+u/3>0$, $y'-v/3-2u/3>0$ and $1-u/3+v/3>0$. We are interested in the zeros of $g_*$ satisfying, additionally, $u>0$, $v>0$ and $z'>0$, corresponding to the original variables $u,v,w,x,y,z$ all being positive. \\\\ 
$g_*$ (unlike $f_*$) is defined when $\eta = 0$ and $g_*(x',y',z',0,0,0)=0$. In Part~\ref{pointfirstbound} we find an upper bound $\eta_1$ on $\eta$ ensuring, via the IFT, that for $x',y'$ and $z'$ lying in the compact positive set termed $\mathcal{Z}$ in the proof, and $0<\eta<\eta_1$, we can solve $g_*=0$ for $u$ and $v$ in terms of $x',y',z'$ and $\eta$ (although $z'$ does not explicitly occur here). In other words, we can locally describe the zero set of $g_*$ as the graph of a function, termed $\theta$ in the proof, with input $(x', y', z', \eta)$ and output $(u,v)$. Moreover, in the region where this description holds $x=x'-v/3+u/3 > 0$, $y=y'-v/3-2u/3>0$, $z=z'>0$ and $w= 1-u/3+v/3 > 0$. The fact that we solve $g_*=0$ not just in the neighbourhood of a single point, but in a neighbourhood of a larger set necessitates use of the slight extension to the IFT in Lemma~\ref{lemIFT1}.\\\\ 
The construction in Part~\ref{pointsecondbound} of the proof gives another upper bound $\eta_2$ on $\eta$ which ensures that $\theta$ is a positive function. In other words, for $0 < \eta < \eta_2$ the solutions to $g_*=0$ obtained in the previous part additionally satisfy $u>0$ and $v>0$. Thus for $x',y'$ and $z'$ lying in the compact positive set $\mathcal{Z}$, and $\eta$ sufficiently small and positive, the graph of $\theta$ corresponds to a positive subset of the original $(x,y,z,u,v,w)$-space.\\\\
To understand the third upper bound on $\eta$ in Part~\ref{pointthirdbound}, consider the ``reduced'' vector field in a portion of $(x',y',z')$-space defined by solving the two algebraic equations of (\ref{Dalg}) for $u$ and $v$ in terms of $x',y',z'$ (possible by the previous arguments for each small enough, fixed, positive $\eta$) and substituting into the differential equations of (\ref{Dalg}). Notice that since $g_*(x',y',z',0,0,0)=0$, the values of $u$ and $v$ on the graph of $\theta$ are small for small $\eta$, and consequently the reduced vector field is close to the original vector field in (\ref{origCRN}). The third upper bound on $\eta$ ensures that the two are close enough to ensure, by regular perturbation theory arguments, that the periodic orbit of (\ref{origCRN}) survives in the reduced system and retains the same number of multipliers inside and outside the unit circle. Thus the reduced system has a periodic orbit of the same stability type as the original vector field (\ref{origCRN}) for sufficiently small values of $\eta$.\\\\
The fourth and final upper bound on $\eta$ in Part~\ref{pointfourthbound} is, roughly speaking, a normal hyperbolicity condition. To understand it, fix any sufficiently small $\eta>0$ and consider (\ref{SingLim}), namely the fast-time system in the limit $\epsilon \to 0$. Note that setting $\dot u = \dot v = 0$ in (\ref{SingLim}) defines the set of equilibria of (\ref{SingLim}), which includes a portion of the graph of the function $\theta$ constructed above. Each equilibrium on this graph has a set of eigenvalues $0$ corresponding to directions tangential to the graph and a set of ``nontrivial'' eigenvalues corresponding to directions transverse to the graph. This final upper bound on $\eta$ ensures that all of the nontrivial eigenvalues have negative real part, and so this portion of the equilibrium manifold of (\ref{SingLim}) is locally, exponentially attracting. This allows use of the perturbation theory results of \cite{Fenichel79}: it ensures, roughly speaking, that the manifold defined by the graph of $\theta$ survives for sufficiently small $\epsilon > 0$ in (\ref{fullCRN}) and has sufficiently smooth dependence on $\epsilon$. Consequently, by regular perturbation theory, the periodic orbit of (\ref{Dalg}) on this graph survives for (\ref{fullCRN}) with sufficiently small $\epsilon > 0$.\\\\
Since the original periodic orbit of (\ref{origCRN}) was linearly stable, so is the perturbed one for (\ref{fullCRN}) -- the previous two bounds on $\eta$ took care of its Floquet multipliers tangential to, and transverse to, the invariant set on which it lies respectively. Thus for sufficiently small, positive $\eta$ and $\epsilon$, (\ref{fullCRN}) has a linearly stable periodic orbit. Moreover, for small enough $\epsilon > 0$, this periodic orbit lies in a region of $(x',y',z',u,v)$-space corresponding to positive $(x,y,z,u,v,w)$-space, since it can be made arbitrarily close to the periodic orbit of the limiting system (\ref{Dalg}). We have thus obtained, as desired, a periodic orbit of (\ref{newCRN}) which is positive and linearly stable relative to its stoichiometry class. 
}

\section{Final remarks}

We remark, first, that giving the added reactions mass action kinetics was convenient and simplified many calculations, but was not fundamental to the techniques of proof of Theorem~\ref{mainthm}. The motivated reader could, with some effort, reprove the result with the added reactions having kinetics from other classes than mass action. Key to the proof is scalability of the reaction rates, and the characterisation of the equilibrium set of the added reactions as a graph over $\mathcal{Z}$. 

The techniques of the proof of Theorem~\ref{mainthm} also provide an alternative proof of Theorem~5 in \cite{banajipanteaMPNE}. The set-up requires only minor and formal modifications: $\mathcal{Z}$ is now a compact subset of some stoichiometry class of $\mathcal{R}$ containing two nondegenerate (resp., linearly stable) equilibria (rather than a periodic orbit), and we need to apply Theorems~12.1~and~12.2 of \cite{Fenichel79} (rather than Theorems~13.1~and~13.2). Using the same approach and the very general and powerful Theorem~9.1 in \cite{Fenichel79}, essentially any compact limit set of (\ref{eqbasic0}), hyperbolic relative to its stoichiometry class, survives for (\ref{perteq}) under the assumptions of Theorem~\ref{mainthm}.

The condition that $\beta$ has rank $m$ can be rephrased, roughly, as ``the new species feature nondegenerately in the new reactions''. Although this condition is essential to our proof \textcolor{black}{(at several points we invert $\hatt{\beta}$), it is possible that a weaker condition might suffice here. The difficulty in exploring this question arises partly because techniques for proving the nonexistence of periodic orbits are limited.}

In Theorem~6 of \cite{banajipanteaMPNE} we proved that nondegenerate (resp., linearly stable) multistationarity is preserved by ``splitting'' reactions and adding in intermediate complexes involving new species, provided this is done in a way satisfying a nondegeneracy condition very similar to the condition on the rank of $\beta$ in Theorem~\ref{mainthm}. We believe an analogous result for periodic orbits to hold, and certain special cases are easily proved. However, the result in full generality cannot be proved using the arguments in the proof of Theorem~\ref{mainthm} here. \textcolor{black}{The set up in} Theorem~6 of \cite{banajipanteaMPNE} does not afford sufficient freedom to control reaction rates to apply the techniques of proof used in Theorem~\ref{mainthm} here: the construction involving two independently controlled parameters, $\eta$ and $\epsilon$ cannot be simply reused, and an alternative approach needs to be found. This task will be undertaken in future work. 

\section*{Acknowledgements}

I would like to thank Casian Pantea, Amlan Banaji, \textcolor{black}{and the anonymous reviewers of this paper} for helpful discussions and useful comments on drafts of this paper.

\end{document}